\documentclass{amsart}
\usepackage{amssymb}
\usepackage{amscd}
\usepackage{verbatim}
\usepackage{epsfig}
\begin{document}

\newcommand\Mand{\ \text{and}\ }
\newcommand\Mwith{\ \text{with}\ }
\newcommand\Mfor{\ \text{for}\ }
\newcommand\Mst{\ \text{such that}\ }
\newcommand\Mor{\ \text{or}\ }
\newcommand\Mif{\ \text{if}\ }
\newcommand\Miff{\ \text{iff}\ }
\newcommand\Mthen{\ \text{then}\ }
\newcommand\nin{\notin}
\newcommand\Id{\operatorname{Id}}
\newcommand\Real{\mathbb{R}}
\newcommand\Cx{\mathbb{C}}
\newcommand\im{\operatorname{Im}}
\newcommand\re{\operatorname{Re}}
\newcommand\Cinf{{\mathcal C}^{\infty}}
\newcommand\dist{{\mathcal C}^{-\infty}}
\newcommand\dCinf{\dot\Cinf}
\newcommand\ddist{\dot\dist}
\newcommand\Sch{{\mathcal S}}
\newcommand\temp{\Sch^{\prime}}
\newcommand\supp{\operatorname{supp}}
\newcommand\Xb{\bar X}
\newcommand\tr{\operatorname{tr}}
\newcommand\Tr{\operatorname{Tr}}
\newcommand\diag{\operatorname{diag}}
\newcommand\vol{\operatorname{vol}}
\newcommand\sing{\text{sing}}
\newcommand\reg{\text{reg}}
\newcommand\Scl{{\text{Sc}}}
\newcommand\scl{{\text{sc}}}
\newcommand\WF{\operatorname{WF}}
\newcommand\WFSc{\WF_\Scl}
\newcommand\WFScp{\WF'_\Scl}
\newcommand\WFsc{\WF_\scl}
\newcommand\pa{\partial}
\newcommand\spp{\operatorname{sp}}
\newcommand\ep{\epsilon}
\newcommand\pih{\hat\pi}
\newcommand\psit{\tilde\psi}
\newcommand\phit{\tilde\phi}
\newcommand\Rn{\Real^n}
\newcommand\sphere{\mathbb{S}}
\newcommand\Sn{\sphere^{n-1}}
\newcommand\Snp{\sphere^n_+}
\newcommand\sci{{}^{\text{sc}}}
\newcommand\sct{\sci T^*}
\newcommand\scdt{\sci \dot T^*}
\newcommand\Fr{{\mathcal F}}
\newcommand\Frinv{\Fr^{-1}}
\newcommand\bop{{\mathcal B}}
\newcommand\spec{\operatorname{spec}}
\newcommand\pspec{\spec_{pp}}
\newcommand\cspec{\spec_{c}}
\newcommand\Psop{\operatorname{\Psi}}
\newcommand\PsiSc{\Psop_\Scl}
\newcommand\PsiSch{\Psop_{\Scl,h}}
\newcommand\PsiScc{\Psop_{\text{scc}}}
\newcommand\calM{{\mathcal M}}
\newcommand\calX{{\mathcal X}}
\newcommand\calC{{\mathcal C}}
\newcommand\SP{\operatorname{RC}}
\newcommand\cl{\operatorname{cl}}
\newcommand\interior{\operatorname{int}}
\newcommand\Span{\operatorname{span}}
\newcommand\scHg{{}^{\text{sc}}H_g}

\setcounter{secnumdepth}{3}
\newtheorem{lemma}{Lemma}[section]
\newtheorem{prop}[lemma]{Proposition}
\newtheorem{thm}[lemma]{Theorem}
\newtheorem{cor}[lemma]{Corollary}
\newtheorem{result}[lemma]{Result}
\newtheorem*{thm*}{Theorem}
\newtheorem*{prop*}{Proposition}
\newtheorem*{conj*}{Conjecture}
\numberwithin{equation}{section}
\theoremstyle{remark}
\newtheorem{rem}[lemma]{Remark}
\theoremstyle{definition}
\newtheorem{Def}[lemma]{Definition}
\newtheorem*{Def*}{Definition}

\renewcommand{\theenumi}{\roman{enumi}}
\renewcommand{\labelenumi}{(\theenumi)}

\title[Spectral shift function in many-body scattering]
{Smoothness and high energy asymptotics
of the spectral shift function in many-body scattering}
\author{Andr\'as Vasy}
\address{Department of Mathematics, Massachusetts Institute of Technology,
Cambridge MA 02139, U.S.A.}
\author{Xue-Ping Wang}
\address{D\'ept.\ Math\'ematiques, Universit\'e de Nantes,
F-44322 Nantes Cedex 03, France}
\date{June 25, 2001}
\thanks{A. V. is partially supported by NSF grant \#DMS-99-70607, and
he thanks the hospitality of the Universit\'e de Nantes, where
this project started. Both authors are grateful
for the hospitality of the Erwin Schr\"odinger
Institute.}
\begin{abstract}
Let $H=\Delta+\sum_{\#a=2} V_a$ be a 3-body Hamiltonian,
$H_a$ the subsystem Hamiltonians, $\Delta\geq 0$ the Laplacian of the
Euclidean metric $g$ on $X_0=\Real^n$, $V_a$ real-valued. Buslaev
and Merkurev \cite{Buslaev-Merkurev:Relationship, Buslaev-Merkurev:Trace}
have shown that, if the pair potentials decay sufficiently fast,
for $\phi\in\Cinf_c(\Real)$, the operator
$\phi(H)-\phi(H_0)-\sum_{\#a=2}(\phi(H_a)-\phi(H_0))$ is trace class.
Hence, one can define a modified
spectral shift
function $\sigma$, as a distribution on $\Real$, by taking its trace.
In this paper we show that if $V_a$ are Schwartz, then
$\sigma$ is in fact $\Cinf$ away
from the thresholds, and obtain its high energy asymptotics.
In addition, we generalize this result to $N$-body scattering, $N$
arbitrary.
\end{abstract}
\maketitle

\section{Introduction}
Let $H=\Delta+\sum_{\#a=2} V_a$ be a 3-body Hamiltonian,
$H_a$ the subsystem Hamiltonians, $\Delta\geq 0$ the Laplacian of the
Euclidean metric $g$ on $X_0=\Real^n$, $V_a$ real-valued. Buslaev
and Merkurev \cite{Buslaev-Merkurev:Relationship, Buslaev-Merkurev:Trace}
have shown that, if the pair potentials decay sufficiently fast,
and under a spectral assumption,
for $\phi\in\Cinf_c(\Real)$, the operator
$\phi(H)-\phi(H_0)-\sum_{\#a=2}(\phi(H_a)-\phi(H_0))$ is trace class.
Indeed, this follows from the adoptation of the
Helffer-Sj\"ostrand functional calculus \cite{Helffer-Sjostrand:Schrodinger},
see \cite{Vasy:Propagation-Many}, even without the spectral assumption.
Hence, one can define a modified
spectral shift function, which is really a distribution on
$\Real_\lambda$, by
\begin{equation}\label{eq:sigma-def}
\sigma(\phi)=\tr(\phi(H)-\phi(H_0)-\sum_{\#a=2}(\phi(H_a)-\phi(H_0))).
\end{equation}
In the two-body setting, $\sigma$ is often denoted by $\xi'$.
The purpose of this paper is to show that $\sigma$ is in fact $\Cinf$ away
from the thresholds, and to obtain its high energy asymptotics. Namely,
we prove the following theorem.

\begin{thm*}(See Theorem~\ref{thm:3-body}.)
Suppose that the pair potentials $V_a$ are Schwartz (on $X^a$).
Then the spectral shift
function $\sigma$ is $\Cinf$ on $\Real\setminus\Lambda$,
where $\Lambda$ is the set of thresholds and $L^2$ eigenvalues of $H$.
Moreover,
$\sigma$ is a symbol outside a compact set, and it has a full asymptotic
expansion as $\lambda\to+\infty$:
$$
\sigma(\lambda)\sim \sum_{j=0}^\infty \lambda^{\frac{n}{2}-3-j} c_j.
$$
In addition,
$$
c_0=C_0\sum_a \sum_{b\neq a}\int_{X_0} V_a V_b\,dg,
$$
where $C_0=\frac{1}{16}(n-2)(n-4)(2\pi)^{-n}\vol(\Sn)$
depends only on $n=\dim X_0$,
and $dg$ is the Riemannian density of the
metric $g$.
\end{thm*}

Buslaev and Merkurev had shown previously \cite{Buslaev-Merkurev:Relationship}
that $\sigma$ is given by a continuous function, under the assumption that
$0$ is not an eigenvalue or resonance of any two-body subsystem, and
that the bottom of the essential spectrum is not an eigenvalue of the
whole Hamiltonian.
We remark that the leading term has an additional $\lambda^{-1}$ decay as
compared to the corresponding formulae in two-body scattering, see
\cite{Buslaev:Trace, Schrader:High, Colin-dV:Formule,
Robert:Asymptotique-EP, Jensen:Phase,
Christiansen:Weyl}.
This is due
to the fact that $\phi(H)-\phi(H_0)-\sum_{\#a=2}(\phi(H_a)-\phi(H_0))$ is
lower order in a high-energy sense than $\phi(H)-\phi(H_0)$ is in the
two-body setting. Correspondingly, if we include a three-particle interaction
$V_1$ in $H$, i.e.\ allow $\#a=1$, the leading term will have its usual
order $\lambda^{\frac{n}{2}-2}$.

In fact,
this is not the only modified spectral shift function one can consider.
For example, let $\chi_a$, $\#a\geq 2$, be a partition of unity by
$\Cinf$ functions on the radial (or geodesic) compactification $\Xb_0$ of the
configuration space $X_0$ that correspond to the collision plane structure.
Then $\phi(H)-\sum_{\#a\geq 2}\phi(H_a)\chi_a$ is also trace class, and
one can define a modified spectral shift function by taking its trace. This
spectral shift
function will depend on the choice of the partition of unity, and
it may be argued that it is less natural than the one adopted above.
This trace is also
smooth away from the thresholds and its high energy asymptotics can also be
calculated. In fact, this expression generalizes without changes to many-body
scattering with arbitrarily many particles, and our proof shows its
smoothness. Even the Buslaev-Merkurev expression can be adopted
to $N$-body scattering, $N$ arbitrary, with combinatorial complexity
being the only additonal issue. The corresponding result is stated
in Theorem~\ref{thm:N-body-trace} in the last section.

The main reason why these statements hold is that the corresponding
expressions for the spectral measure, or for the high-energy cutoff outgoing
(or incoming)
resolvents $\psi(H)(H-(\lambda+i0))^{-1}$, $\psi\in\Cinf_c(\Real)$,
already make sense pointwise in $\lambda$. While these operators
are not necessarily trace class, the trace makes sense as an oscillatory
integral since it can be regarded as a pairing with the delta
distribution associated to the diagonal.
Indeed, the latter only requires that the scattering wave front
set of the kernel of the the operator whose trace we intend to take, is
disjoint from the conormal bundle of the diagonal lifted to the b-double
space $[\Xb_0\times\Xb_0;\pa\Xb_0\times\pa\Xb_0]$. In view of the propagation
of singularities in many-body scattering, applied to the resolvent kernel, 
this is automatically satisfied provided that the kernel of the operator
combination that we consider is essentially a generalized eigenfunction
of some many-body Hamiltonian microlocally near the conormal bundle of the
diagonal.

Our methods are thus an adaptation of traditional microlocal analysis
to many-body scattering. In fact, we do not need the whole program
initiated by Melrose \cite{RBMSpec}, see \cite{Vasy:Bound-States},
since the trace is a rather simple object. Thus, a moderate strengthening
of propagation results in
the dual of the radial variable essentially suffices.
These have their origins in
the Mourre estimate \cite{Mourre-Absence}, and became partially
microlocal estimates in the work of G\'erard, Isozaki, Skibsted
\cite{GerComm, GIS:N-body} and the work of the second author \cite{XPWang};
see also \cite{Hassell-Vasy:Symbolic} for a discussion of these in the
geometric setting.

The weak high energy asymptotics, i.e.\ that of
$$
\tr(\phi(H/\lambda)-\phi(H_0/\lambda)
-\sum_{\#a=2}(\phi(H_a/\lambda)-\phi(H_0/\lambda))),
$$
$\phi\in\Cinf_c(\Real)$ fixed, can be derived from the semiclassical
functional calculus of Helffer-Robert
\cite{Robert:Calcul, Helffer-Robert:Calcul}. Symbol estimates for
$\sigma$ follow from combining microlocal versions of
high energy estimates for
$R(\lambda+i0)$, \cite{Jensen:High},
with our method of proving smoothness of $\sigma$ (see also \cite{Wang:High}).
The perturbation series expansion for $R(\lambda+i0)$ then gives the
full asymptotics.

The assumption that the $V_a$ are Schwartz is not optimal. The
distributional trace, $\sigma$, is defined if the $V_a$ are symbols of
order $<-n$ (on $X^a$). Below we also obtain partial results for these
potentials and explain the finer tools one needs to extend the Theorem
to this setting. Specifically, this would entail obtaining, at least
microlocally near the conormal bundle of the diagonal at infinity,
oscillatory integral estimates in place of the wave front set estimates.

We also mention that under stronger assumptions as in
\cite{Gerard:Distortion}, such as analyticity of
the pair potentials in a cone near infinity, one can use complex scaling,
or distortion analyticity, to continue $(H-\lambda)^{-1}$ across the
spectrum. Since the trace is invariant under the scaling, and the
scaled operators are elliptic near the real axis, we deduce that the
trace is real analytic.

While there is quite a bit of machinery in the microlocal
analysis of this paper, due to the complexities of many-body scattering,
given the microlocal estimates,
the arguments proving the Theorem are essentially the same as in two-body
scattering. Hence we invite the reader who is familiar with the
scattering microlocal analysis introduced by Melrose \cite{RBMSpec}
to read Section~\ref{sec:smooth} pretending that one works in the
two-body setting. Moreover, in the appendix we sketch the proof of
a propagation theorem that is weaker than the results of
\cite{Vasy:Bound-States}, but suffices for the purposes of this paper.
The sketch is intended to stand on its own, and to explain the basic
ideas in microlocalizing propagation estimates.

The structure of the paper is as follows. In Section~\ref{sec:notation}
we recall the standard many-body notation and make some preliminary
remarks. In the following section we remind the reader of the microlocal
structure of many-body scattering. In Section~\ref{sec:smooth} we prove
the smoothness of $\sigma$ in the three-body setting,
and in Section~\ref{sec:high} we obtain the high
energy expansion of $\sigma$. In the final section we state and prove the
corresponding results in many-body scattering, and in the appendix we
outline the proof of the simplified positive commutator estimate
described above.

A.\ V.\ is grateful to Maciej Zworski for numerous stimulating
discussions on traces
in the geometric two-body setting (as arising in the work of Christiansen
\cite{Christiansen:Weyl}). The authors also thank
Leonid Friedlander, Richard Melrose, Didier Robert
and Steve Zelditch for helpful comments.

\section{Notation and preliminaries}\label{sec:notation}
Before we can state the
precise definitions, we need to introduce some basic (and mostly
standard) notation. We refer to \cite{Derezinski-Gerard:Scattering} for
a very detailed discussion of the setup and the basic results.
We consider the Euclidean space $\Rn$, and let
$g$ be the standard Euclidean metric
on it. We assume also that we
are given a (finite) family $\calX$ of linear subspaces
$X_a$, $a\in I$,
of $\Rn$ which is closed under intersections and includes the subspace
$X_1=\{0\}$ consisting of the origin, and the whole space $X_0=\Rn$.
Let $X^a$ be the
orthocomplement of $X_a$.
We write $g_a$ and $g^a$ for the induced metrics on $X_a$ and $X^a$
respectively.
We let $\pi^a$ be the orthogonal
projection to $X^a$, $\pi_a$ to $X_a$. A many-body Hamiltonian is an operator
of the form
\begin{equation}
H=\Delta+\sum_{a\in I} (\pi^a)^*V_a;
\end{equation}
here $\Delta$ is the positive Laplacian, $V_0=0$,
and the $V_a$ are real-valued
functions in an appropriate class.

There is a natural partial
ordering on $I$ induced by the ordering of $X^a$
by inclusion, i.e.\ $b\leq a$ means that $X^b\subset X^a$, hence
$X_b\supset X_a$.
A three-body Hamiltonian is a many-body Hamiltonian with
$I\neq\{0,1\}$ such that for all $a,b\nin\{0,1\}$
with $a\neq b$, $X_a\cap X_b=\{0\}$ holds.

Corresponding to each cluster $a$ we introduce the cluster
Hamiltonian $H^a$ as an operator on $L^2(X^a)$ given by
\begin{equation}
H^a=\Delta_{X^a}+\sum_{b\leq a} V_b,
\end{equation}
$\Delta_{X^a}$ being the Laplacian of the induced metric on $X^a$.
With
\begin{equation}
\calX^a=\{X_b\cap X^a:\ b\leq a\},
\end{equation}
$H^a$ is a many-body Hamiltonian with collision planes $\calX^a$.
We write
\begin{equation}
H_a=H^a\otimes\Id_{X_a}+\Id_{X^a}\otimes \Delta_{X_a}
=\Delta_{X_0}+\sum_{b\leq a} V_a.
\end{equation}
The $L^2$ eigenfunctions of $H^a$, also called bound states,
can be used to define the set of thresholds of $H^b$. Namely, we let
\begin{equation}
\Lambda_a=\cup_{b<a}\pspec(H^b)
\end{equation}
be the set of thresholds of $H^a$, and we also let
\begin{equation}
\Lambda'_a=\Lambda_a\cup\pspec(H^a)
=\cup_{b\leq a}\pspec(H^b).
\end{equation}
Thus, $0\in\Lambda_a$ for $a\neq 0$ and $\Lambda_a\subset(-\infty,0]$.
It follows from the Mourre
theory (see e.g.\ \cite{FroMourre, Perry-Sigal-Simon:Spectral})
that $\Lambda_a$ is compact,
countable, and $\pspec(H^a)$ can only accumulate at $\Lambda
_a$. We also let
\begin{equation}
\Lambda=\Lambda_1\cup\pspec(H).
\end{equation}
For future reference we also define the intercluster interactions as
\begin{equation}
I_a=H-H_a=\sum_{b\not\leq a}V_b.
\end{equation}

One of the properties of many-body operators that we exploit below is
that if $H_L$ and $H_R$ are many-body Hamiltonians on $(X_0)_L,(X_0)_R$
respectively, then
\begin{equation}
P=P_\alpha=\alpha H_L\otimes\Id+(1-\alpha)\Id\otimes H_R,\qquad 0<\alpha<1,
\end{equation}
is a many-body
Hamiltonian on $M_0=(X_0)_L\times(X_0)_R$. (We name the new space
$M_0$ and the operator $P$ to avoid confusing notation.)
The collision planes are
given by $(X_a)_L\times (X_b)_R$, $a,b\in I$, though they have a rather
special structure generated by $(X_a)_L\times (X_0)_R$ and
$(X_0)_L\times (X_a)_R$, $a\in I$, by taking intersections of these.
The thresholds $\Lambda=\Lambda_\alpha$
of $P$ are of the form $\alpha\lambda_L+(1-\alpha)\lambda_R$,
where $\lambda_L$ and $\lambda_R$ are thresholds of $H_L$ and $H_R$
respectively.

Since we want to study the Schwartz kernel of the resolvent
on a space, for us the relevant case is $X_0=(X_0)_L=(X_0)_R$
(this is the space),
and $H_L=H_R$, so the kernel of $(H_L-\lambda)^{-1}$ is a distribution on
$M_0=(X_0)_L\times(X_0)_R$. The thresholds of $P_\alpha$ vary with $\alpha$,
except that $\lambda\in\Lambda_R=\Lambda_L$ implies that
$\lambda\in\Lambda_\alpha$ for all $\alpha$. On the other hand, given
$\lambda\nin\Lambda_L$, there is a countable subset $C$ of $(0,1)$ such
that $\alpha\nin C$ implies $\lambda\nin\Lambda_\alpha$. Indeed, for any
pair $(\lambda_L,\lambda_R)\in\Lambda_L\times\Lambda_R$
with $\lambda_L\neq\lambda_R$, there is
a unique $\alpha\in\Real$ such that
$\lambda=\alpha\lambda_L+(1-\alpha)\lambda_R$, and then we can take $C$ to be
the intersection of the set of such $\alpha$'s with $(0,1)$.

In the next section we explain the full microlocal propagation picture
proved in \cite{Vasy:Bound-States}. However, we would like to emphasize
here that apart from mild technical issues we only need an estimate
of G\'erard, Isozaki, Skibsted \cite{GerComm, GIS:N-body} and the
second author \cite{XPWang}, and its slightly strengthened version.
Namely, let $B=\frac{w}{2\langle w\rangle}\cdot D_w
+ D_w\cdot \frac{w}{2\langle w\rangle}$,
and suppose that $F_1,F_2\in\Cinf_c(\Real)$, $\supp F_1\subset (c,c')$,
$\supp F_2\subset (c',c'')$, $c<c'<c''$, $|c|,|c''|$ sufficiently small,
and $\psi\in\Cinf_c(\Real)$.
The aforementioned papers show that
$F_1(B) \psi(H)(H-(\lambda+i0))^{-1} F_2(B)$ maps
$H^{0,s}$ to $H^{0,s'}$ for all $s,s'$
when $\lambda\nin\Lambda$. In the terminology of the next
section, this holds because the principal symbol of $B$ at infinity, which
is denoted by $-\tau$ there, is increasing along generalized broken
bicharacteristics, and singularities propagate in the forward direction
under $(H-(\lambda+i0))^{-1}$.

We only need the following
strengthening of this statement (as well as its high energy version).
Suppose that $\chi_1,\chi_2\in S^0(X_0)$ are conic cutoff functions (outside
a compact subset of $X_0$) with $\supp\chi_1\cap\supp\chi_2=\emptyset$.
Then there exists $\ep>0$ (depending on the supports and on
$\lambda\nin\Lambda$) such that for
$F_1,F_2\in\Cinf_c(\Real)$, $\supp F_1\subset (c,c'+\ep)$,
$\supp F_2\subset (c',c'')$, $c<c'<c''$, $|c|,|c''|$ sufficiently small,
and $\psi\in\Cinf_c(\Real)$,
$\chi_1 F_1(B) \psi(H)(H-(\lambda+i0))^{-1} F_2(B)\chi_2$
maps $H^{0,s}$ to $H^{0,s'}$ for all $s,s'$
when $\lambda\nin\Lambda$. In the microlocal
terminology, this means that the principal symbol of $B$ at infinity
has to increase by at least $\ep>0$ along generalized broken bicharacteristics
leaving $\supp\chi_2$ (near infinity) and ending in $\supp\chi_1$.
We remark that the smoothness order $r$ of the Sobolev space $H^{r,s}$
is left as $0$ only due to the non-ellipticity (in the usual sense) of $B$.
Slightly different pseudo-differential cutoffs, such as those described
in the next section, yield maps from $H^{r,s}$ to $H^{r',s'}$ for all
$r,r',s,s'$.

We also need the high energy version of both of these estimates, which
can be proved similarly. Namely, first, that for all $r,s,r',s'$,
\begin{equation}
\|F_1(B) \psi(H)(H-(\lambda+i0))^{-1} F_2(B)\|
_{\bop(H^{0,s},H^{0,s'})}\leq C\lambda^{-1/2},\quad\lambda\geq 1,
\end{equation}
if $F_1,F_2\in\Cinf_c(\Real)$, $\supp F_1\subset (c,c')$,
$\supp F_2\subset (c',c'')$, $c<c'<c''$, $|c|,|c''|$ sufficiently small,
and $\psi\in\Cinf_c(\Real)$. Next,
there exists $\ep>0$ (depending on the supports of $\chi_1$ and
$\chi_2$ as above), but not on
$\lambda\geq 1$) such that for
$F_1,F_2\in\Cinf_c(\Real)$, $\supp F_1\subset (c,c'+\ep)$,
$\supp F_2\subset (c',c'')$, $c<c'<c''$, $|c|,|c''|$ sufficiently small,
and $\psi\in\Cinf_c(\Real)$, and
\begin{equation}
\|\chi_1 F_1(B) \psi(H)(H-(\lambda+i0))^{-1} F_2(B)\chi_2\|
_{\bop(H^{0,s},H^{0,s'})}\leq C\lambda^{-1/2},\quad\lambda\geq 1.
\end{equation}

\section{Microlocal analysis in many-body scattering}
We recall the microlocal tools introduced in the many-body setting in
\cite{Vasy:Propagation-Many}. Thus,
we compactify $X_0$ as in \cite{RBMSpec} by letting
\begin{equation}
\Xb=\Xb_0=\Snp
\end{equation}
to be the radial compactification of $X_0$ (also called the
geodesic compactification) to a closed
hemisphere, i.e.\ a ball. Identifying $X_0$ with $\Rn$ (and using $g$)
we write the compactification map $\SP:\Rn\to\Snp$ given by
\begin{equation}\label{eq:SP-def}
\SP(w)=(1/(1+|w|^2)^{1/2},w/(1+|w|^2)^{1/2})\in\Snp\subset\Real^{n+1},
\quad w\in\Rn.
\end{equation}
One can view the use of the compactification as using inverted
polar coordinates $(r^{-1},\omega)\in[0,1)\times\Sn$
near infinity, i.e.\ working near `$r^{-1}=0$', where $w=r\omega$
in usual polar coordinates.
We write the coordinates on $\Rn=X_a\oplus X^a$ as $(w_a,w^a)$.
We let
\begin{equation}
\Xb_a=\cl(\SP(X_a)),\quad C_a=\Xb_a\cap\partial\Snp.
\end{equation}
Hence, $C_a$ is a sphere of dimension $n_a-1$ where
$n_a=\dim X_a$. We also let
\begin{equation}
\calC=\{C_a:\ a\in I\}.
\end{equation}
Thus, $C_0=\partial\Snp=\Sn$,
and $a\leq b$ if and only if $C_b\subset C_a$. We also define the regular
part of $C_a$ as
\begin{equation}
C_{a,\reg}=C_a\setminus\cup_{b\not\leq a} C_b=C_a\setminus \bigcup\{C_b:
\ C_b\subsetneq C_a\}
\end{equation}
Since throughout this paper we work in the Euclidean setting,
where the notation $X$, $X_a$, etc., has been used for the (non-compact)
vector spaces, we always use a bar, as in $\Xb$, $\Xb_a$, etc., to
denote the corresponding compact spaces. We write $H^r$ for the standard
Sobolev space, and $H^{r,s}$ for the weighted Sobolev space
$\langle w\rangle^{-s}H^r$.

In \cite{Vasy:Propagation-Many} a pseudo-differential operator calculus
was constructed on many-body spaces as above $(\Rn,\calX)$;
indeed, it was defined
for appropriate geometric generalizations of these spaces. Operators
of multi-order $(k,l)$ are denoted by $\PsiSc^{k,l}(\Xb_0;\calC)$;
the notation only refers to $\calC$ rather than $\calX$ since this is
the only pertinent information in the geometric setting.
We briefly recall the definition of the many-body pseudo-differential
calculus via the quantization of symbols. Below $[\Xb_0;\calC]$ is the
blow-up of $\Xb_0$ at $\calC$; in particular, its interior is diffeomorphic
to that of $\Xb_0$, i.e.\ to $X_0=\Rn$, and there is a smooth blow-down map
$[\Xb_0;\calC]\to\Xb_0$, so every $\Cinf$ function on $\Xb_0$ is $\Cinf$
on $[\Xb_0;\calC]$.
Thus, the polyhomogeneous space
$\PsiSc^{m,l}(\Xb_0,\calC)$ is the following.
We identify $\interior(\Xb_0)$ and $\interior([\Xb_0;\calC])$ with $\Rn$
as usual (via $\SP^{-1}$), suppose that
$a\in\Cinf(\Rn_w\times\Rn_\xi)$ is in fact of the form
\begin{equation}\label{eq:cl-symb-def}
a\in\rho_{\infty}^{-m}\rho_{\partial}^l\Cinf([\Xb_0;\calC]\times\Xb_0^*),
\end{equation}
where $\rho_\partial$ and $\rho_\infty$ are defining functions of the first
and second factors, $\Xb_0$ and $\Xb_0^*$, respectively, so they can be
taken as $\langle w\rangle^{-1}$ and $\langle \xi\rangle^{-1}$ respectively.
Let $A=q_L(a)$ denote the
left quantization of $a$:
\begin{equation}\label{eq:q_L-def}
Au(w)=(2\pi)^{-n}\int e^{i(w-w')\cdot\xi}a(w,\xi)u(w')\,dw'\,d\xi,
\end{equation}
understood as an oscillatory integral. Then $A\in\PsiSc^{m,l}(\Xb_0,\calC)$.
We could have equally well used other (right, Weyl, etc.) quantizations as
well, and we could have also allowed $a$ to depend on $w'$ as well
(with $w'$ regarded as a smooth variable on $[\Xb_0;\calC]$).

An advantage of using symbols
\begin{equation}\label{eq:cl-symb-double}
a\in\Cinf([\Xb_0;\calC]_w\times[\Xb_0;\calC]_{w'} \times(\Xb_0^*)_\xi)
\end{equation}
is that certain operators with amplitudes $a=q(w,\xi)p(w',\xi)$, which
can be considered the operator product $A=QP$
of the left quantization $Q$ of $q$
and the right quantization $P$ of $p$, lie in the
class $\PsiSc^{m,l}(\Xb_0,\calC)$ even though the left quantization of
$q$ does not. As an example, for $\psi_0\in\Cinf_c(\Real)$,
we can write $\psi_0(H)$ as the right quantization of a symbol $p$ as above
that is in fact Schwartz in $\xi$.
Then $q$ does not have to lie $\Cinf([\Xb_0;\calC]\times\Xb_0^*)$ to make
\eqref{eq:cl-symb-double} hold, since
due to the rapid decay of $p$ in $\xi$, the $\xi\to\infty$ behavior
of $q$ is (mostly) irrelevant.
For example, we can allow $q=f(w\cdot\xi/|w|)$,
$f\in\Cinf_c(\Real)$; then
$a=qp$ satisfies \eqref{eq:cl-symb-double}. Thus, $A=A(f)$ defined by
\begin{equation}\label{eq:q_W-def}
Au(w)=(2\pi)^{-n}\int e^{i(w-w')\cdot\xi}a(w,w',\xi)u(w')\,dw'\,d\xi,
\quad a(w,w',\xi)=q(w,\xi)p(w',\xi),
\end{equation}
yields an
element of $\PsiSc^{-\infty,l}(\Xb_0,\calC)$. The operator $A(f)$ has
very similar properties to $f(B)\psi_0(H)$, $B$ as in the previous section,
and can be considered as a pseudo-differential replacement for $f(B)\psi_0(H)$,
for $f(B)\psi_0(H)$ is not a ps.d.o.\ (due to the lack of ellipticity of $B$).
(Note though that it follows from the arguments of \cite{Hassell-Vasy:Symbolic}
that $\psi_1(H)f(B)\psi_2(H)
\in\PsiSc^{-\infty,0}(\Xb_0,\calC)$ for $\psi_1,\psi_2\in\Cinf_c(\Real)$.)

Elements of $\PsiSc^{k,l}(\Xb_0;\calC)$
are in particular bounded operators from $H^{r,s}$ to $H^{r-k,s+l}$.
(The different signs are due to the index convention introduced in the
geometric scattering setting in \cite{RBMSpec}.) It was also shown
in \cite{Vasy:Propagation-Many}
that for $\lambda\nin\spec(H)$, $R(\lambda)=(H-\lambda)^{-1}
\in\PsiSc^{-2,0}(\Xb_0,\calC)$, and the Helffer-Sj\"ostrand argument
was adopted to show that for $\phi\in\Cinf_c(\Real)$, $\phi(H)\in
\PsiSc^{-\infty,0}(\Xb_0,\calC)$, so it is smoothing in the standard sense,
but does not give any decay at infinity (see Section~\ref{sec:high} for
a sketch of this argument in the semiclassical setting).
The construction of $R(\lambda)$
is local at $\pa\Xb_0$ (in fact, it is microlocal in a sense discussed below).
Thus, if $H=\Delta+V$, $H'=\Delta+V'$ are many-body operators such that in
a neighborhood of $p\in\pa\Xb_0$, $H-H'=V-V'$ is in $x^k\Cinf(\Xb_0)$,
$k\geq 1$, i.e.\ $\chi(V-V')\in x^k\Cinf(\Xb_0)$ for some
$\chi\in\Cinf(\Xb_0)$, $\chi(p)\neq 0$,
then $R_H(\lambda)-R_{H'}(\lambda)$ is in
$\PsiSc^{-2,k}(\Xb_0,\calC)$ near $p$, so $(R_H(\lambda)-R_{H'}(\lambda))\chi
\in\PsiSc^{-2,k}(\Xb_0,\calC)$. This in turn implies, via the
Helffer-Sj\"ostrand construction, that for $\phi\in\Cinf_c(\Real)$,
$\phi(H)-\phi(H')$ is in $\PsiSc^{-\infty,k}(\Xb_0,\calC)$ near $p$.

Applying this in our setting, i.e.\ for 3-body Hamiltonians $H$ with
potential $\sum_{\#a=2}V_a$, note that if the potentials $V_b$ are
classical symbols of order $-k$ on $X^b$, then $H-H_a$ is in
$x^k\Cinf(\Xb_0)$ away from $\cup_{\#b=2,\ b\neq a} C_b$, and
$H_a-H_0$ is in $x^k\Cinf(\Xb_0)$ away from $C_a$. Hence
$\phi(H)-\phi(H_a)$ is in $\PsiSc^{-\infty,k}(\Xb_0,\calC)$ away from
$\cup_{\#b=2,\ b\neq a} C_b$, and $\phi(H_a)-\phi(H_0)$ is in
the same class away from $C_a$. Now consider the expression
$\phi(H)-\phi(H_0)-\sum_{\#a=2}(\phi(H_a)-\phi(H_0))$ near $C_b$. Rewriting
it as $\phi(H)-\phi(H_b)-\sum_{\#a=2,\ a\neq b}(\phi(H_a)-\phi(H_0))$
shows that it is in $\PsiSc^{-\infty,k}(\Xb_0,\calC)$ away from
$\cup_{\#a=2,\ a\neq b} C_a$, hence in particular near $C_b$.
Since $b$ is arbitrary, we deduce that
$\phi(H)-\phi(H_0)-\sum_{\#a=2}(\phi(H_a)-\phi(H_0))\in
\PsiSc^{-\infty,k}(\Xb_0,\calC)$. Since pseudo-differential operators
of order $(m,l)$, $m<-n$, $l>n$, are easily seen to be of trace class
(the kernel is in particular continuous and integrable along the diagonal),
we deduce that $\phi(H)-\phi(H_0)-\sum_{\#a=2}(\phi(H_a)-\phi(H_0))$ is
trace class
indeed, as discussed by Buslaev and Merkurev \cite{Buslaev-Merkurev:Trace}.
This allows us to define $\sigma$ by taking its trace, i.e.\ by
\eqref{eq:sigma-def}.

We need more refined tools to analyze the smoothness of $\sigma$. Namely,
we need to analyze the Schwartz kernels of operators that do not lie
in $\PsiSc^{m,l}(\Xb_0,\calC)$, e.g.\ of $(H-(\lambda+i0))^{-1}$.
This requires the use of some
microlocal analysis, which we briefly recall below.

The phase space in scattering theory is the cotangent bundle
$T^*X_0$ which can be identified with $X_0\times(X_0)^*$.
Again, following Melrose \cite{RBMSpec},
it is convenient to consider its appropriate partial
compactification, i.e.\ to consider it as a trivial vector bundle over $\Xb_0$
by compactifying the base. Hence, one defines the scattering cotangent
bundle of $\Xb_0$ by $\sct\Xb_0=\Xb_0\times X_0^*$.
We remark that the construction of $\sct\Xb_0$
is completely natural and geometric, just like the following ones,
see \cite{RBMSpec}. Note that $(X_0)^*$ can be identified with $X_0$
via the metric $g$ as we often do below.

Recall from \cite{Vasy:Propagation-Many} that
in many-body scattering $\sct\Xb_0$ is
{\em not} the natural place for microlocal analysis for the very same
reason that introduces the compressed cotangent bundle
in the study of the wave equation on bounded domains; see the
works of Melrose, Sj\"ostrand
\cite{Melrose-Sjostrand:I} and Lebeau \cite{Lebeau:Propagation} on
the wave equation in domains with smooth boundaries or corners,
respectively.
We can see what causes trouble from both the dynamical and
the quantum point of view. Regarding dynamics,
the issue is that only the external part of the momentum
is preserved in a collision, the internal part is not; while from
the quantum point of view the problem is that there is only partial
commutativity in the algebra of
the associated pseudo-differential operators, even to top
order. Hence, one cannot expect to localize in arbitrary open subsets of
$\sct_{\pa\Xb_0}\Xb_0$, i.e.\ to microlocalize fully.
To rectify this, we replace the full bundle
$\sct_{C_{a,\reg}}\Xb_0=C_{a,\reg}\times(X_0)^*\simeq
C_{a,\reg}\times X_0$
over $C_{a,\reg}\subset\pa\Xb_0$ 
by $\sct_{C_{a,\reg}}\Xb_a=C_{a,\reg}\times X_a$, i.e.\ we consider
\begin{equation}
\scdt \Xb_0=\cup_a \sct_{C_{a,\reg}}\Xb_a.
\end{equation}
Over $C_{a,\reg}$, there is a natural projection $\pi_a:\sct_{C_{a,\reg}}\Xb_0
\to\sct_{C_{a,\reg}}\Xb_a$, i.e.\ $\pi_a:C_{a,\reg}\times X_0^*
\to C_{a,\reg}\times X_a^*$, corresponding to the pull-back of one-forms;
in the trivialization given by the metric it is
induced by the orthogonal projection to $X_a$ in
the fibers. By putting the $\pi_a$ together, we obtain a projection
$\pi:\sct_{\pa\Xb_0}\Xb_0\to\scdt\Xb_0$. We put the topology induced by
$\pi$ on $\scdt\Xb_0$. As mentioned above,
this definition is analogous to that of the compressed
cotangent bundle in the study of
the wave equation in domains with smooth boundaries or corners.

We also recall from \cite{RBMSpec} that the characteristic variety
$\Sigma_0(\lambda)$
of $\Delta-\lambda$ is simply the subset of $\sct_{\pa\Xb_0}\Xb_0$ where
$g-\lambda$ vanishes; $g$ being the metric function. If $\Lambda_1=\{0\}$,
the compressed characteristic set of $H-\lambda$ is
$\pi(\Sigma_0(\lambda))\subset\scdt\Xb_0$. In general, all the bound
states contribute to the characteristic variety. Thus, we let
\begin{equation}
\Sigma_b(\lambda)=\{\zeta=(y_b,\xi_b)\in\sct_{C_b}\Xb_b:
\ \lambda-|\xi_b|^2\in\pspec{H^b}\}\subset\sct_{C_b}\Xb_b;
\end{equation}
note that $|\xi_b|^2$ is the kinetic energy of a particle in a bound
state of $H^b$.
If $C_a\subset C_b$, there is also a natural projection
$\pi_{ba}:\sct_{C_{a,\reg}}\Xb_b\to\sct_{C_{a,\reg}}\Xb_a$
(in the metric trivialization we can use the orthogonal projection
$X_b\to X_a$ as above), and then we define
the characteristic set of $H-\lambda$ to be
\begin{equation}
\dot\Sigma(\lambda)=\cup_a\dot\Sigma_a(\lambda),\qquad\dot\Sigma_a(\lambda)
=\cup_{C_b\supset C_a}
\pi_{ba}(\Sigma_b(\lambda))\cap\sct_{C_{a,\reg}}\Xb_a,
\end{equation}
so $\dot\Sigma(\lambda)\subset\scdt \Xb_0$.

The radial subset $R_+(\lambda)\cup R_-(\lambda)$ of $\dot\Sigma(\lambda)$
plays a special role in many-body scattering, for points in this
form stationary generalized broken bicharacteristics. Namely,
that at points $(w/|w|,\xi)\in R_\pm(\lambda)$,
the Hamilton vector field of $\Delta_{X_b}$ is radial (so there is no
propagation tangentially to $T^* X_b$ as discussed below), and simultaneously
the particles may be in a bound state of $H^b$, hence there is
no propagation in normal directions either.
More explicitly, let
$-\tau$ be the dual of the radial Euclidean variable, so
\begin{equation}
\tau=-\frac{w\cdot\xi}{|w|}.
\end{equation}
Then
\begin{equation}\begin{split}
R_\pm(\lambda)=\{\xi&=(\omega,\xi_a)\in\sct_{C_{a,\reg}}\Xb_a:\\
&\ \exists b,\ C_a\subset C_b,
\ \lambda-\tau(\xi)^2\in\pspec(H^b),\ \pm\tau(\xi)=|\xi_a|\}
\end{split}\end{equation}
are the incoming ($+$) and outgoing ($-$) radial sets respectively.
Note that if $\lambda$ is not a threshold of $H$, then $\tau\neq 0$ on
$R_+(\lambda)\cup R_-(\lambda)$. The set $R_+(\lambda)\cup R_-(\lambda)$
is the propagation set of Sigal and Soffer \cite{Sigal-Soffer:N};
this is the set where there is no real principal type propagation.

Generalized broken bicharacteristics of a many-body Hamiltonian $H-\lambda$
were defined in
\cite[Definition~2.1]{Vasy:Bound-States}.
Here we do not recall the full definition here, but remind the reader
that these are continuous maps $\gamma:I\to\dot\Sigma(\lambda)$, where
$I\subset\Real$ is an interval, that enjoy certain estimates with
regard to Hamilton vector fields in the various subsystems. The main
property for us is that the function
$\tau=-\frac{w\cdot\xi}{|w|}$ on $\scdt \Xb_0$,
is strictly decreasing along generalized broken bicharacteristics except
at the radial points. Note that $-\tau$ is the dual of the radial variable,
and as mentioned in the introduction, it played a major role in previous
works on many-body scattering.

As mentioned in the introduction, `singularities' (i.e.\ lack of
decay at infinity) of $u\in\temp$ are described
by the many-body scattering wave front set, $\WFSc(u)$, which was
introduced in \cite{Vasy:Propagation-Many}, and which describes $u$ modulo
Schwartz functions, similarly to how the usual wave front set
describes distributions modulo smooth functions. Just as for the image
of the bicharacteristics, $\scdt\Xb_0$ provides the natural setting
in which $\WFSc$ is defined: $\WFSc(u)$ is a closed subset of $\scdt\Xb_0$.
The definition of $\WFSc(u)$ relies on the algebra of many-body scattering
pseudo-differential operators, $\PsiSc^{k,l}(\Xb_0,\calC)$.
There are several possible definitions of $\WFSc$,
all of which agree for generalized eigenfunctions of $H$, but the
one given in \cite{Vasy:Propagation-Many} that is modelled on the
fibred-cusp wave front set of Mazzeo and Melrose \cite{Mazzeo-Melrose:Fibred}
enjoys many properties of the usual wave front set.

We remark that in
the two-body setting, when $\scdt\Xb_0=\sct_{\pa\Xb_0}
\Xb_0$, $\WFSc$ is just the
scattering wave front set $\WFsc$ introduced by Melrose, \cite{RBMSpec},
which in turn is closely related to the usual wave front set via the Fourier
transform. Thus, for $(\omega,\xi)\in\sct_{\pa\Xb_0}\Xb_0$, considered as
$\pa\Xb_0\times\Rn=\partial\Xb_0\times\Rn$,
$(\omega,\xi)\nin\WFsc(u)$ means that there exists
$\phi\in\Cinf(\Xb_0)$ such that $\phi(\omega)\neq 0$ and $\Fr(\phi u)$ is
$\Cinf$ near $\xi$. If we employed the usual conic terminology instead
of the compactified one, we would think of $\phi$ as a conic cut-off
function in the direction $\omega$. Thus, $\WFsc$ at infinity is
analogous to $\WF$ with the role of position and momentum reversed.
The definition of $\WFSc(u)$ is more complicated, but if
$u=\psi(H)v$ for some $\psi\in\Cinf_c(\Real)$ (any other operator
in $\PsiSc^{-\infty,0}(\Xb,\calC)$ would do instead of $\psi(H)$), then
the following is
a sufficient condition for $(\omega,\xi_a)\in\sct_{C_{a,\reg}}\Xb_a$,
considered as $C_{a,\reg}\times X_a$, not
to be in $\WFSc(u)$. Suppose that there exists
$\phi\in\Cinf(\Xb_0)$, $\phi(\omega)\neq 0$, and $\rho\in\Cinf_c(X_a)$,
$\rho(\xi_a)\neq 0$, and $((\pi_a)^*\rho)\Fr(\phi u)\in\Sch(\Rn)=\Sch(X_0)$.
Then $(\omega,\xi_a)\nin\WFSc(u)$.
This wave front set also has complete analogues
for the relative wave front set $\WFSc^{k,l}(u)$, i.e.\ when one is
working modulo weighted Sobolev spaces $H^{k,l}$. In addition, although
the definition of $\scdt\Xb_0$ depends on the metric $g$, there is
a natural identification of the compressed spaces corresponding to different
metrics $g$, and the wave front set is defined independently of the choice
of $g$ then.

As an example, consider $M_0=X_0\times X_0$ as the space, and let $u\in
\temp(M_0)$ be the Schwartz kernel of $\psi(H)$, $H$ a
many-body Hamiltonian on $X_0$, $\psi\in\Cinf_c(\Real)$. As discussed
above $M_0$ is naturally a many-body space with collision planes
$X_a\times X_b$, $a,b\in I$. However, to accommodate the diagonal
singularity, we introduce additional collision planes by adding
the diagonal $\diag=\{(w,w):\ w\in X_0\}$ to this set, as well as
the diagonals of $X_a$: $\diag_a=\diag\cap(X_a\times X_a)$,
$a\in I$; the resulting set is denoted by $\calM$. Recall that
$\overline{\diag}$ is the closure of $\diag$ in $\bar M_0$, hence
it is the b-diagonal as discussed in \cite{RBMSpec, Vasy:Propagation-Many}.
Now the kernel $K$
of $\psi(H)$ is conormal to $\pa\overline{\diag_a}$
in the strong sense that it becomes
$\Cinf$ upon the blow-up of these (in the appropriate order); indeed,
this is essentially the definition of $\PsiSc^{-\infty,0}(\Xb_0;\calC)$.
Note that if $\psi_0\in\Cinf_c(\Real)$ is identically $1$ on
$[\inf\supp\psi,\sup\supp\psi]$, then $\psi_0(P_\alpha)K=K$, as follows
directly from the functional calculus since $H\otimes\Id$, $\Id\otimes H$
commute, so the Fourier transform description of the wave front set
given above is applicable.

Since $K$ is the kernel of a pseudo-differential operator, it is immediate
that $\WFSc(K)\subset\cup_a \sct_{\pa\overline{\diag}_{a,\reg}}
\overline{\diag_a}$, i.e.\ it lies over the boundary of the diagonal; indeed,
the kernel is rapidly decreasing elsewhere by the definition
of $\PsiSc^{-\infty,0}(\Xb_0;\calC)$.
We claim that for all $a$,
$\WFSc(K)\cap\sct_{\pa\overline{\diag}_{a,\reg}}
\overline{\diag_a}$ lies in the zero section of
$\sct_{\pa\overline{\diag}_{a,\reg}}\overline{\diag_a}$. To see this,
let $\phi\in\Cinf(\bar M_0)$ be such that $\supp\phi$ is disjoint
from all collision planes other than those containing $\diag_a$.
This means that $\phi$ localizes near the regular, and away from the
singular, part of $\pa\overline{\diag}_a$. The orthocomplement
of $\diag_a$, with respect to $g_L+g_R$,
is $\diag_{X^a\times X^a}\oplus \{(w,-w):\ w\in X_0\}$, and
the corresponding orthogonal projections give coordinates
$(w_a+w'_a, w^a+(w')^a,w-w')$ (up to a factor of $\sqrt{2}$).
The corresponding dual variables are (up to another factor)
$(\xi_a+\xi'_a, \xi^a+(\xi')^a,\xi-\xi')$.
Since $\phi K$ is Schwartz in $w-w'$,
and is symbolic in $w+w'$ after some blow ups, in particular symbolic
in $w_a+w'_a$, the Fourier transform of $\phi K$ is Schwartz outside
$\xi_a+\xi'_a=0$. Hence $\WFSc(K)$ lies in $\xi_a+\xi'_a=0$, i.e.\ in
the zero section of
$\sct_{\pa\overline{\diag}_{a,\reg}}\overline{\diag_a}$ as claimed.

A similar argument proves that $\WFSc(\delta_{\diag})$ lies in
the zero section of
$\sct_{\pa\overline{\diag}_{a,\reg}}\overline{\diag_a}$. Indeed,
one can see this directly from the definition, for any vector field
tangent to $\diag$ annihilates $\delta_{\diag}$.

Note that the orthogonal projection of the conormal bundle of $\diag$
to $T^*\diag_a$ is the zero section, and for this reason we call
the zero section of
$\sct_{\pa\overline{\diag}_{a,\reg}}\overline{\diag_a}$ the
compressed conormal bundle of $\diag$ at $\diag_a$. We summarize these
statements as a lemma.

\begin{lemma}\label{lemma:psdo-kernel-WFSc}
The kernel $K$ of $\psi(H)$, $\psi\in\Cinf_c(\Real)$, satisfies
$\WFSc(K)$ is a subset of the zero section of
$\sct_{\pa\overline{\diag}_{a,\reg}}\overline{\diag_a}$, i.e.\ of the
compressed conormal bundle of the diagonal at $\diag_a$. The same holds for
$\delta_{\diag}$. In particular, $\tau=0$ on $\WFSc(K)$ and on
$\WFSc(\delta_{\diag})$.
\end{lemma}

\begin{rem}
For example, consider $A(F)\in\PsiSc^{-\infty,0}(\bar M_0,\pa\calM)$ defined
by \eqref{eq:q_W-def}, $F\in\Cinf_c(\Real)$,
and note that on its operator wave front
set (a notion that we briefly recall in the Appendix; see \cite[Section~5]
{Vasy:Propagation-Many} for a detailed discussion) $-\tau\in\supp F$.
The lemma immediately implies that for $F\in\Cinf_c(\Real)$ with
$0\nin\supp F$,
$A(F)K\in\Sch$.
\end{rem}

The main property of wave front sets that we use is the following.
If $u,v\in\temp$, $\WFSc(u)$ is a compact subset
of $\scdt\Xb_0$, and $\WFSc(u)\cap\WFSc(v)=\emptyset$, then the $L^2$
pairing extends by continuity to define $\langle u,v\rangle$.
For the relative wave front set one has stronger results. In particular,
it suffices that $\WFSc^{*,l}(u)\cap\WFSc^{*,-l}(v)=\emptyset$; the
smoothness order $k$ is irrelevant since $\WFSc(u)$ is assumed to be compact.
This property follows immediately from
\cite[Proposition~5.4]{Vasy:Propagation-Many} and the definition of $\WFSc$.

The theorem on the propagation of singularities is the following.

\begin{thm*}(\cite[Theorem~2.2]{Vasy:Bound-States})
Let $u\in\temp(\Rn)$, $\lambda\in\Real$. Then
\begin{equation}
\WFSc(u)\setminus\WFSc((H-\lambda)u)
\end{equation}
is a union of maximally extended generalized
broken bicharacteristics of $H-\lambda$ in
$\dot\Sigma(\lambda)\setminus\WFSc((H-\lambda)u)$.
\end{thm*}

A stronger version of this is valid for $u=(H-(\lambda+i0))^{-1}f$,
defined if $\WFSc(f)$ satisfies appropriate
non-incoming conditions (see \cite[Section~2]{Vasy:Bound-States} for
a full discussion). 
Namely $\WFSc(f)$ only propagates forward along generalized broken
bicharacteristics. The non-incoming conditions
are satisfied, for example, if $\tau\leq 0$
on $\WFSc(f)$, and then we conclude that $\tau\leq 0$ on $\WFSc(u)$.
These statements are also proved by positive commutator estimates, and
consequences can be easily restated as operator estimates.
For example, consider $A(F)\in\PsiSc^{-\infty,0}(\Xb_0,\calC)$ defined
by \eqref{eq:q_W-def}.
The fact that $-\tau$ is increasing
along generalized broken bicharacteristics then converts into the estimate
of \cite{GerComm, XPWang}. Indeed, mapping properties of the
wave front set, \cite[Proposition~5.4]{Vasy:Propagation-Many}, imply
that for $\lambda$ in a compact subset of $\Real\setminus\Lambda$,
$r,s,r',s'$ arbitrary,
\begin{equation}\label{eq:WF-B-8}
\|A(F) R(\lambda+i0)A(\tilde F)\|_{\bop(H^{r,s},H^{r',s'})}\leq C,
\end{equation}
where $\tilde F\in\Cinf_c(\Real)$, and in
addition $\supp F\subset (c,c')$, $\supp \tilde F\subset (c',c'')$,
$c<c'<c''$. $|c|, |c'|$ sufficiently small.
The fact that in addition $-\tau$ is strictly increasing on
any non-stationary generalized broken bicharacteristic implies the
strengthened estimate that for all $\chi,\tilde\chi\in\Cinf(\Xb_0)$ with
$\supp\chi\cap\supp\tilde\chi\cap\pa\Xb_0=\emptyset$, and
for $\lambda$ in a compact subset of $\Real\setminus\Lambda$,
$r,s,r',s'$ arbitrary, there exists $\ep>0$ such that
\begin{equation}\label{eq:WF-B-16}
\|\chi A(F) R(\lambda+i0)
A(\tilde F)\tilde\chi\|
_{\bop(H^{r,s},H^{r',s'})}\leq C,
\end{equation}
where $\tilde F\in\Cinf_c(\Real)$, and in
addition $\supp F\subset (c,c'+\ep)$, $\supp \tilde F\subset (c',c'')$,
$c<c'<c''$. $|c|, |c'|$ sufficiently small. These are very similar to
the estimates mentioned at the end of Section~\ref{sec:notation}.
Of course, these estimates can be simplified, but here we explicitly
wanted to emphasize the connection between $\WFSc$ and the estimates.
We outline the proof of \eqref{eq:WF-B-16} in the appendix, so that the
reader can avoid the more complicated arguments leading to the
stronger results of \cite{Vasy:Bound-States}.

\section{Smoothness of the trace}\label{sec:smooth}
Let $K=K_\lambda,K_a=(K_a)_\lambda$ 
denote the kernels of $R(\lambda+i0)=(H-(\lambda+i0))^{-1}$, and
$R_a(\lambda+i0)=(H_a-(\lambda+i0))^{-1}$ respectively, in $\im\lambda\geq 0$
(of course $+i0$ can be dropped if $\im\lambda>0$). Hence they are
tempered distributions on $\Xb_0\times\Xb_0$. Then $K$ solves
$(H_L-\lambda)K=\delta_{\diag}$, $(H_R-\lambda)K=\delta_{\diag}$,
where $H_L$, $H_R$ act on the left and right factors of $X_0\times X_0$
respectively, and $\delta_{\diag}$ is the kernel of the identity
operator, hence is the natural
Dirac delta distribution associated to the diagonal.
(The second equation really contains the transpose of
$H_R-\lambda$, but since $H_R$ is self-adjoint and real, this is just $H_R$.)
Thus, for example,
\begin{equation}
(\alpha H_L+(1-\alpha)H_R-\lambda)K_\lambda=\delta_{\diag},\quad 0<\alpha< 1,
\end{equation}
for $\im\lambda\geq 0$. Hence, for $\im\lambda>0$,
\begin{equation}
K_\lambda=(\alpha H_L+(1-\alpha)H_R-\lambda)^{-1}\delta_{\diag}.
\end{equation}
As remarked in Section~\ref{sec:notation},
$P_\alpha=\alpha H_L+(1-\alpha)H_R$
is a many-body Hamiltonian, so all the usual results on
many-body Hamiltonians, such as propagation of singularities, are applicable.
In particular, if $\lambda\nin\Lambda$, we can arrange that
$\lambda\nin\Lambda_\alpha$
by an appropriate choice of $\alpha$.
We let $\tau=-(w\cdot\xi+w'\cdot\xi')/\langle(w,w')\rangle$ as usual, with
$(w,w')$ denoting coordinates on $X_0\times X_0$, and $(\xi,\xi')$ the
dual coordinates. Thus, $-\tau$ is the dual of the radial variable
$r=|(w,w')|$ with respect to the metric $g_L+g_R$. We thus deduce that
$(\alpha H_L+(1-\alpha)H_R-(\lambda+i0))^{-1}\delta_{\diag}$
makes sense since $\tau=0$ on $\WFSc(\delta_{\diag})$ by
Lemma~\ref{lemma:psdo-kernel-WFSc}. Thus,
\begin{equation}
K_\lambda=(\alpha H_L+(1-\alpha)H_R-(\lambda+i0))^{-1}\delta_{\diag}
\equiv (P_\alpha-(\lambda+i0))^{-1}\delta_{\diag}.
\end{equation}

For technical reasons it is convenient to eliminate the diagonal singularity
(which cancels from the spectral measure anyway). This can be done by
considering
\begin{equation}
\psi(H)(H-(\lambda+i0))^{-1}=\psit(H)
\end{equation}
instead of $R(\lambda+i0)$,
where $\psi\in\Cinf_c(\Real)$ is
identically $1$ near $\lambda$. This has similar properties to $R(\lambda+i0)$,
except that it does not have a diagonal singularity in the interior.
In particular, its kernel, which we denote by $\tilde K=\tilde K_\lambda$,
satisfies $(H_L-\lambda)\tilde K=(H_R-\lambda)\tilde K=K_{\psi(H)}$,
$K_{\psi(H)}$ denoting the kernel of $\psi(H)$. Hence, again using
Lemma~\ref{lemma:psdo-kernel-WFSc},
\begin{equation}
\tilde K_\lambda=(\alpha H_L+(1-\alpha)H_R-(\lambda+i0))^{-1}K_{\psi(H)}
\equiv (P_\alpha-(\lambda+i0))^{-1}K_{\psi(H)}.\end{equation}

As recalled in the previous section,
$\tau$ is decreasing under the forward generalized
broken bicharacteristic relation.
Then, by the propagation of singularities, namely
that they only propagate in the forward direction under
$(P_\alpha-(\lambda+i0))^{-1}$, we deduce that
$\tau\leq 0$ on $\WFSc(\tilde K)$. Similar statements hold
for $\tilde K_a$. Hence $\tau\leq 0$ on $\WFSc(\tilde K')$,
\begin{equation}
\tilde K'=\tilde K-\tilde K_0-\sum_a (\tilde K_a-\tilde K_0).
\end{equation}
We can write this explicitly by saying that if $F\in\Cinf_c(\Real)$,
supported in $(-\infty,0)$, then $A(F)\tilde K\in\Sch$, where now
$A(F)$ is an operator on $M_0=X_0\times X_0$, with principal symbol
at infinity the `same' as that of $F(B)\tilde\psi(H)$ (though the latter
is not quite a ps.d.o.), $B$
the radial vector field on $X_0\times X_0$, with principal symbol $-\tau$.
Moreover, $\tau$ is strictly negative everywhere
on this wave front set except possibly at the compressed conormal bundle of the
diagonal, i.e.\ at the zero section of
$\sct_{\pa\overline{\diag}_{a,\reg}}\overline{\diag_a}$.

We remark that while microlocally on $\WFSc(\tilde K)$, $\tilde K$ is not
Schwartz, it is of course
still in a weighted Sobolev space, namely in
$H^{r,s-1}$ for any $s<-n/2$ and for any $r$. This holds because
$K_{\psi(H)}$ is in $H^{r,s}$ for any $s<-n/2$ and for any $r$, and
$(P_\alpha-(\lambda+i0))^{-1}$ loses one order of decay.
(This is the standard propagation phenomenon as for the wave equation, and
the factor $h^{-1}$ appears in \eqref{eq:WF-B-16-h} for the same reason.)

We claim
that in fact $\tau\leq\tau_0< 0$ everywhere, hence in particular the
compressed conormal bundle of the diagonal is disjoint from it.
The key step is the following lemma.

\begin{lemma}
Let $b$ be a 2-cluster. Then $\WFSc(\tilde K-\tilde K_b)$ and
$\WFSc(\tilde K_a-\tilde K_0)$, $a\neq b$,
are disjoint from the zero section of
$\sct_{\pa\overline{\diag}_{b,\reg}}\overline{\diag_b}$ and of
$\sct_{\pa\overline{\diag}_{0,\reg}}\overline{\diag_0}$.
\end{lemma}

\begin{proof}
To see this, we first compute
$(P_\alpha-\lambda)u$, $u=\tilde K-\tilde K_b$.
Thus,
\begin{equation}\label{eq:prop-sing-32}
(P_\alpha-\lambda)u=K_{\psi(H)}-K_{\psi(H_b)}
+[\alpha (I_b)_L+(1-\alpha)(I_b)_R] \tilde K_b.
\end{equation}
Thus, near $\pa\overline{\diag_b}$,
the result is rapidly decreasing since the intercluster interactions
$(I_b)_L$, $(I_b)_R$ are
such, and the same holds for
the difference $K_{\psi(H)}-K_{\psi(H_b)}$. {\em It is only here that we use
$V_a\in \Sch(X^a)$}; see the comments of the following paragraphs for
potentials that do not decay rapidly. Hence, in this region,
$\WFSc(u)$ is a union of
maximally extended generalized broken bicharacteristics of
$P_\alpha-\lambda$. That is, if $\zeta$ is in
$\sct_{\pa\overline{\diag}_{b,\reg}}\overline{\diag_b}\cup
\sct_{\pa\overline{\diag}_{0,\reg}}\overline{\diag_0}$, and
$\zeta\in\WFsc(u)$, then there exists $T>0$ and a generalized
broken bicharacteristic $\gamma:[-T,T]\to\dot\Sigma(\lambda)$ such that
$\gamma(0)=\zeta$ and $\gamma(t)\in\WFSc(u)$ for all $t\in[-T,T]$.
(The constant $T$ depends on $\zeta$, namely on how long it takes
for bicharacteristics to leave
$\sct_{\pa\overline{\diag}_{b,\reg}}\overline{\diag_b}\cup
\sct_{\pa\overline{\diag}_{0,\reg}}\overline{\diag_0}$.) But if
$\zeta$ is in addition in the zero section of
$\sct_{\pa\overline{\diag}_{b,\reg}}\overline{\diag_b}$ or of
$\sct_{\pa\overline{\diag}_{0,\reg}}\overline{\diag_0}$, then
$\tau(\gamma(0))=0$, and $\tau(\gamma(t))>0$ for $t<0$, so
$\gamma(t)\nin\WFSc(u)$ (as $\tau\leq 0$ on $\WFSc(u)$) contradicting
that $\gamma(t)\in\WFSc(u)$ for all $t\in[-T,T]$.
Here we use again
that $\lambda$ is not a threshold, for this assumption ensures that there are
no stationary bicharacteristics at $\tau=0$.
Thus, the zero section of
$\sct_{\pa\overline{\diag}_{b,\reg}}\overline{\diag_b}$ and of
$\sct_{\pa\overline{\diag}_{0,\reg}}\overline{\diag_0}$, i.e.\ 
the compressed
conormal bundle of the diagonal, is disjoint from $\WFSc(u)$.
A similar argument shows
that for $b\neq a$, $\WFSc(\tilde K_a-\tilde K_0)$
is disjoint from the compressed conormal bundle
of the diagonal in the same region.
\end{proof}

Combining the statements of the lemma shows
that $\WFSc(\tilde K')$ has the same property for every $b$. Since $b$ is
arbitrary, this means that $\WFSc(\tilde K')$ is (globally) disjoint
from the compressed conormal bundle of the diagonal, and $\tau$ is negative
on the former. (The existence of a strictly negative constant $\tau_0$
such that $\tau\leq\tau_0$ on $\WFSc(\tilde K')$ follows by compactness.)
In explicit terms, this means that there exists $\ep>0$ such that
if $F\in\Cinf_c(\Real)$,
supported in $(-\infty,\ep)$, then $A(F)\tilde K'\in\Sch$.

The trace of $\tilde K'$, defined
as the push-forward of the restriction of $\tilde K'$ to the diagonal,
thus makes
sense by wave-front set considerations. Note that the trace is thus just
the $L^2$-pairing of $\tilde K'$ with $\delta_{\diag}$.
Hence the only
property of $\WFSc$ we need is that distributions with disjoint $\WFSc$
can be $L^2$-paired, which we already discussed in the previous section.
Again, explicitly this argument amounts to writing
\begin{equation}\label{eq:explicit-pairing}
\langle\tilde K',\delta_{\diag}\rangle
=\langle (\Id-\psi_0(P_\alpha))
\tilde K',\delta_{\diag}\rangle+
\langle A(F)\tilde K',\delta_{\diag}\rangle
+\langle \tilde K',(\psi_0(P_\alpha)-A(F)^*)\delta_{\diag}\rangle,
\end{equation}
$F$ as in the previous paragraph, with in addition $F\equiv 1$ near $0$,
and noting that both terms on the right hand side are defined. (Rather than
using continuity from $\Sch$ in the style of H\"ormander \cite{FIO1},
one could use continuity from $\im\lambda>0$ here.)

Note that this argument also shows that the trace of
$\chi_b (\psi(H)R(\lambda+i0)-\psi(H_b)R_b(\lambda+i0))$
is already defined as an oscillatory
integral, or rather by wave front considerations.

Applying the same arguments for $R(\lambda-i0)$, we deduce the same pointwise
for the corresponding expressions involving the spectral measure, i.e.\ for
\begin{equation}
\spp(\lambda)-\spp_0(\lambda)-\sum_a(\spp_a(\lambda)-\spp_0(\lambda)),
\end{equation}
where
\begin{equation}
\spp(\lambda)=(2\pi i)^{-1}((H-(\lambda+i0))^{-1}-(H-(\lambda-i0))^{-1}),
\end{equation}
and similarly with the other spectral measures. Since $\phi(H)=
\int_\Real\phi(\lambda)\,\spp(\lambda)\,d\lambda$, we conclude that
the distribution $\sigma$ is given by the continuous function
\begin{equation}
\sigma(\lambda)
=\tr(\spp(\lambda)-\spp_0(\lambda)-\sum_a(\spp_a(\lambda)-\spp_0(\lambda))).
\end{equation}
The smoothness of this function follows directly by taking derivatives
of the resolvents with respect to $\lambda$. Indeed the derivative
of $\tilde K$ with respect to $\lambda$ is given by
\begin{equation}
\frac{d}{d\lambda}
\tilde K_\lambda=(P_\alpha-(\lambda+i0))^{-2}K_{\psi(H)}
=(P_\alpha-(\lambda+i0))^{-1}\tilde K_\lambda.
\end{equation}
Then $\tau\leq 0$ still holds on $\WFSc(\frac{d}{d\lambda}\tilde K_\lambda)$, 
$\WFSc(\tilde K_\lambda)$. In addition,
\begin{equation}
(P_\alpha-\lambda)\frac{d}{d\lambda}\tilde K_\lambda
=\tilde K_\lambda
\end{equation}
shows that
\begin{equation}
(P_\alpha-\lambda)\frac{d}{d\lambda}(\tilde K_\lambda
-(\tilde K_b)_\lambda)=\tilde K_\lambda-(\tilde K_b)_\lambda
+[\alpha (I_b)_L+(1-\alpha)(I_b)_R] (\tilde K_b)_\lambda.
\end{equation}
Since near $\pa\overline{\diag_b}$,
$\WFSc(\tilde K_\lambda-(\tilde K_b)_\lambda)$ is in
$\tau\leq\tau_0<0$, we deduce as above that
$\WFSc(\frac{d}{d\lambda}(\tilde K_\lambda
-(\tilde K_b)_\lambda))$ is also in $\tau\leq\tau_0<0$
near $\pa\overline{\diag_b}$, and then further
that $\WFSc(\frac{d}{d\lambda}\tilde K')$ is in $\tau\leq\tau_0<0$ everywhere.
This shows that $\frac{d\sigma}{d\lambda}$ is continuous, and an iterative
argument yields that $\sigma$ is $\Cinf$.
We have thus proved the following proposition, which forms the first
part of the Theorem stated in the introduction.

\begin{prop}
Suppose that $H$ is a three-body Hamiltonian, and
the pair potentials $V_a$ are Schwartz (on $X^a$).
Then the spectral shift function $\sigma$, defined by \eqref{eq:sigma-def},
is $\Cinf$ on $\Real\setminus\Lambda$.
\end{prop}

If the potentials $V_a$ are not Schwartz on $X^a$, rather they are symbols
in $S^{-\rho}(X^a)$, the first part of \eqref{eq:WF-B-16-h},
namely $K_{\psi(H)}-K_{\psi(H_b)}$, is in $H^{r,s+\rho}$ for any
$s<-n/2$ (and any $r$)
near $\pa\overline{\diag}_b$,
while $[\alpha(I_b)_L+(1-\alpha)(I_b)_R] \tilde K_b$
is in $H^{r,s-1+\rho}$ in the same region.
This yields that $u=\tilde K-\tilde K_b$ is in $H^{r,s-2+\rho}$ near the
compressed conormal bundle of $\overline{\diag}_b$ (rather than Schwartz).
The pairing with $\delta_{\diag}$ thus makes sense if $(s-2+\rho)+s\geq 0$,
i.e.\ $\rho\geq -2s+2$. Since $s>-n/2$ is arbitrary, this means $\rho>n+2$
suffices. This result is non-optimal, $\rho>n$ should be sufficient,
but the improvement should correspond
to more delicate (oscillatory integral type) behavior
that we do not consider here since proving the better estimates
is more complicated, especially if the set of thresholds is large!
In addition, each derivative causes an
additional order of decay to be lost, so to conclude that $\sigma$ is
$\calC^k$ by this argument we need $\rho>n+2+k$.

We explain briefly in terms of the two-body setting why more precise,
oscillatory integral type, control is needed on the resolvent kernel
near the conormal bundle of the diagonal to obtain optimal results.
In two-body scattering, the high-energy cut-off resolvent kernel
$\tilde K$ is of the form $e^{i\sqrt{\lambda}\langle w-w'\rangle}
\langle w-w'\rangle^{-(n-1)/2}a$,
$a$ smooth on the blown up space $[\bar M_0;\pa\overline{\diag}]$
microlocally near the conormal bundle of the diagonal, and is
also a smooth function of $\lambda>0$. Differentiating the kernel,
in particular the exponential, with respect to $\lambda$, gives an extra
factor of growth $\langle w-w'\rangle$. However, this factor is bounded
near the diagonal, so the pairing with $\delta_{\diag}$ makes sense
without imposing extra vanishing conditions on $V$. This argument shows
that $\rho>n$ suffices in the two-body setting to conclude that
$\sigma$ is $\Cinf$; indeed, it immediately
generalizes to the geometric setting
of \cite{RBMSpec}, where it was proved by the use of trace formulae in
\cite{Christiansen:Weyl}, since one also knows
the similar oscillatory behavior there \cite{Hassell-Vasy:Spectral}.
(This FIO-type analysis originated with \cite{RBMZw} in the geometric
setting.)
Such a detailed analysis is much harder in the many-body setting.

Returning to the setting of Schwartz pair interactions,
we remark that a perhaps better approach, which we do not pursue here,
for proving smoothness,
is to include $\lambda$ as a variable, i.e.\ to work with, say, the
Schr\"odinger equation directly. Then the result of the pull-back to
the diagonal in the spatial variables followed by the push-forward
can be seen to be smooth by the wave front set calculus.

We also emphasize that apart from technical issues, the whole argument
can be rewritten in terms estimates involving operators such as $F(B)$,
$B$ as in the previous sections.

\section{High energy asymptotics}\label{sec:high}
Below we change high energy problems into semiclassical problems, so
we start by recalling some notation.
First, $H_h^{r,s}=H_h^{r,s}(\Rn)$ is the semiclassical
weighted Sobolev space, $\langle w\rangle^{-s}H_h^r(\Rn)$, where the norm
on $H_h^r$ is defined by
\begin{equation}
\|u\|_{H_h^r}=(2\pi h)^{-n/2}
\|\langle\xi\rangle^{r}{\mathcal{F}}_h u\|_{L^2(\Rn)},
\end{equation}
${\mathcal{F}}_h$ being the $h$-Fourier transform,
\begin{equation}
({\mathcal{F}}_h u)(\xi)=\int e^{-iw\cdot\xi/h}u(w)\,dw.
\end{equation}
For $r=0$, the norm in $H_h^r(\Rn)$ is the
standard $L^2(\Rn)$ norm; the norm on $H_h^1(\Rn)$ is
$\|(1+h^2\Delta)u\|_{L^2(\Rn)}$, etc.
Semiclassical many-body pseudo-differential operators
$A=A(h)\in\PsiSch^{m,l}(\Xb_0,\calC)$ are defined by
the following analogue of \eqref{eq:q_L-def}:
\begin{equation}\label{eq:semicl-quant}
(2\pi h)^{-n}\int e^{i(w-w')\cdot\xi/h} a(w,\xi;h)\,d\xi
\end{equation}
where $a\in\rho_\infty^{-m}\rho_\pa^l
\Cinf([\Xb_0;\calC]\times X_0^*\times[0,1))$,
rapidly decreasing in the second factor. One can similarly modify
\eqref{eq:q_W-def}, to define $A=A(F;h)$ by
\begin{equation}\label{eq:q_W-def-h}
A(F;h)u(w)=(2\pi h)^{-n}\int e^{i(w-w')\cdot\xi/h}a(w,w',\xi;h)
u(w')\,dw'\,d\xi,
\end{equation}
$a(w,w',\xi;h)=q(w,\xi;h)p(w',\xi;h)$,
this yields an element of $\PsiSch^{-\infty,l}(\Xb_0,\calC)$.
One can also interpret these spaces as high energy spaces by letting
$h=\lambda^{-1/2}$; then $\lambda\to\infty$ corresponds to $h\to 0$.

As an example, we compute the norm of kernel $K_{A(h)}$ of
$A=A(h)
\in\PsiSch^{-\infty,0}(\Xb_0,\calC)$ in $H^{0,s}_h((X_0)_w\times (X_0)_{w'})$
for $s<-n/2$. Namely, \eqref{eq:semicl-quant} takes the form
\begin{equation}
h^{-n} A(w,(w-w')/h;h),\ A\in\Cinf([\Xb_0;\calC]\times X_0\times
[0,1)),
\end{equation}
$A$ rapidly decreasing in the second factor. Thus, for $s<-n/2$
\begin{equation}\begin{split}\label{eq:A(h)-ker-norm}
\|K_{A(h))}\|_{H^{0,s}_h}^2
&=h^{-2n}\int_{X_0\times X_0}\langle(w,w')\rangle^{2s}
|A(w,(w-w')/h;h)|^2\,dw\,dw'\\
&=h^{-n}\int_{X_0\times X_0}\langle w\rangle^{2s}
|A(w,W;h)|^2\,dw\,dW\leq C'h^{-n},
\end{split}\end{equation}
hence $\|K_{A(h)}\|_{H^{0,s}_h}\leq Ch^{-n/2}$. In fact, a similar
calculation, after inserting $(\Id+h^2\Delta_{(w,w')})^r$ in front
of $A$, gives that $\|K_{A(h)}\|_{H^{r,s}_h}\leq Ch^{-n/2}$ for all $r$
(and for all $s<-n/2$).
A similar calculation also applies to $\delta_{\diag}$, and gives
for $s<-n/2$ that
\begin{equation}
\|\delta_{\diag}\|_{H^{s,s}_h}\leq Ch^{-n/2}.
\end{equation}

We start our study of high energy asymptotics of $\sigma$
by stating an averaged version.

\begin{prop}
Suppose that $H$ is a three-body Hamiltonian, $V_a$ are symbols
of order $-k$ on $X^a$, $k>n$. Then for $\phi\in\Cinf_c(\Real)$,
\begin{equation}
\tr(\phi(H/\lambda)-\phi(H_0/\lambda)
-\sum_a (\phi(H_a/\lambda)-\phi(H_0/\lambda)))
\sim \sum_{j=0}^\infty
\lambda^{\frac{n}{2}-2-j} d_j,\qquad\lambda\to +\infty.
\end{equation}
\end{prop}

\begin{proof}
We change the notation slightly and compute
\begin{equation}
\tr(\phi(H/\lambda')-\phi(H_0/\lambda')
-\sum_a (\phi(H_a/\lambda')-\phi(H_0/\lambda')))
\end{equation}
as $\lambda'\to+\infty$.
As usual, we convert the high energy problem into a semiclassical one
by letting $h^2=(\lambda')^{-1}$, and let
\begin{equation}
H(h)=h^2 H,
\end{equation}
so $\phi(H/\lambda')=\phi(H(h))$, i.e.\ we are interested in the
intersection of the spectrum of $H(h)$ with a compact interval (namely
$\supp\phi$). Note that the high energy asymptotics $\lambda'\to+\infty$
corresponds to the semiclassical problem $h\to 0$.
Note also that
$h^2 H=h^2\Delta+h^2 V$, so semiclassically the potential vanishes
two orders higher than the semiclassical Laplacian $h^2\Delta$.
Let $\phit\in\Cinf_c(\Cx)$ be an almost analytic extension of $\phi$
(so $|\bar\partial_{\lambda}\phit(\lambda)|\leq C_N|\im\lambda|^{N}$ for all
$N$),
and let $R(\lambda;h)=(H(h)-\lambda)^{-1}$.
Via the Cauchy formula,
\begin{equation}\label{eq:def-phi(H)}
\phi(H(h))=-\frac{1}{2\pi i}\int\bar\partial_{\lambda}\phit(\lambda)
R(\lambda;h)\,d\lambda\wedge d\bar \lambda,
\end{equation}
the study of the functional calculus in a semiclassical setting
reduces to that of the behavior of the resolvents
away from the real axis, and the uniformity in $h$, up to the real axis
in $\lambda$.
The main fact that makes
the calculations uniform in $h$ is that by self-adjointness,
the $L^2$ operator norm of $R(\lambda;h)$ is bounded by $|\im \lambda|^{-1}$,
for this implies (via the parametrix construction) that all seminorms
of $R(\lambda;h)$ in the semiclassical version of
$\PsiSc^{-2,0}(\Xb_0,\calC)$ (cf.\ \eqref{eq:semicl-quant}) are bounded by
$C|\im\lambda|^{-k}$
for $\re\lambda$ in a compact set, $|\im\lambda|$ bounded from above,
and for some $C$ and $k$ depending on the seminorm.
Hence it follows immediately that for
$\phi\in\Cinf_c(\Real)$,
\begin{equation}\begin{split}\label{eq:fc-tr}
&\tr(\phi(H(h))-\phi(H_0(h))-\sum_a (\phi(H_a(h))-\phi(H_0(h))))\\
&=\tr\left(-\frac{1}{2\pi i}\int\bar\partial_{\lambda}\phit(\lambda)
(R(\lambda;h)-R_0(\lambda;h)-\sum_a (R_a(\lambda;h)-R_0(\lambda;h)))
\,d\lambda\wedge d\bar \lambda\right)\\
&\qquad\sim \sum_{j=0}^\infty h^{-n+j} a_j.
\end{split}\end{equation}

While the integral is trace class, this may not be apparent since
the integrand, which is in $\PsiSc^{-2,k}(\Xb_0,\calC)$,
is not such due to the
diagonal singularity of the resolvents. However, it is easy to rewrite the
integral so that the integrand becomes trace class as well.
Namely, let $\psi(\lambda)=(\lambda-\lambda_0)^{m}\phi(\lambda)$,
$m>(\dim X_0-2)/2$, $\lambda_0\nin\Real$, so
$\phi(H(h))=R(\lambda_0;h)^m\psi(H(h))$, and write out the Cauchy
integral representation of $\psi(H(h))$ using an almost analytic
continuation $\tilde\psi$ of $\psi$. Thus,
\begin{equation}\begin{split}\label{eq:fc-tr-reg}
&\phi(H(h))-\phi(H_0(h))-\sum_a (\phi(H_a(h))-\phi(H_0(h)))\\
&=-\frac{1}{2\pi i}\int\bar\partial_{\lambda}\tilde\psi(\lambda)
(R(\lambda_0;h)^m R(\lambda;h)-R_0(\lambda_0;h)^m R_0(\lambda;h)\\
&\qquad\qquad
-\sum_a (R_a(\lambda_0;h)^m R_a(\lambda;h)-R_0(\lambda_0;h)^m R_0(\lambda;h)))
\,d\lambda\wedge d\bar \lambda.
\end{split}\end{equation}
The integrand is now clearly trace class, in fact it is a semiclassical
many-body scattering pseudo-differential operator of order
$(-2m-2,\rho)$, and the asymptotics as $h\to 0$
follows from the construction of $R(\lambda;h)$, etc.

Our first remark is that since $H$ is a differential operator (rather than
a more general pseudo-differential operator),
for all odd $j$, $a_j=0$, as in \cite{Robert:Asymptotique-EP}.
This follows from the fact that these terms are given
by integrals of odd functions of the momentum $\xi$,
hence these integrals vanish.
Indeed, due to the additional $h^2$ vanishing in $V$, which implies
similar vanishing for the difference of the various resolvents,
the leading term is at least two orders higher, i.e.\ $a_0=0$ (and we
have already mentioned that $a_1=0$).
This is exactly the same order of vanishing as for two-body
scattering. In fact, more is true.

To see this, we perform local calculations near $C_b$, i.e.\ we
replace the full trace by the trace of $\chi_b$ times \eqref{eq:fc-tr-reg}
where $\chi_b\in\Cinf(\Xb_0)$ is supported away from the $C_a$ for $a\neq b$.
Making sure that the $\chi_b$ form a partition of unity, \eqref{eq:fc-tr}
becomes the sum of the local traces, so it suffices to show that
the $a_2$ term of each local trace vanishes.
Thus,
$\chi_b(R(\lambda_0;h)^m R(\lambda;h)-R_b(\lambda_0;h)^m R_b(\lambda;h))$,
$\chi_b(R_a(\lambda_0;h)^m R_a(\lambda;h)-R_0(\lambda_0;h)^m R_0(\lambda;h))$,
$a\neq b$, are already trace class. Using $R(\lambda;h)-R_b(\lambda;h)
=R(\lambda;h)h^2 I_b R_b(\lambda;h)$, a similar formula for $R_a(\lambda;h)
-R_0(\lambda;h)$, that $R(\lambda;h)$, $R_b(\lambda;h)$
can be replaced by $R_0(\lambda;h)$ up to an error of $h^2$, and the various
operators can be commuted up to an error with an additional power of $h$,
it is easy to see that the $a_2$ term vanishes, where we also take into
account that $I_b=\sum_{a\neq b}V_a$. (In addition, $a_3=0$ since the odd
coefficients vanish.)
Changing back to the $h$-independent notation
finishes the proof.

An alternative, but somewhat formal calculation
(in the sense that it does not use
the regularization procedure) proceeds as follows.
Note first that the resolvents can be rewritten as
\begin{equation}\begin{split}
&R(\lambda;h)=R_0(\lambda;h)-\sum_a R_a(\lambda;h) h^2V_a R_0(\lambda;h)\\
&\qquad\qquad\qquad+\sum_a\sum_{b\neq a} R(\lambda;h)h^2 V_b
R_a(\lambda;h)h^2 V_a R_0(\lambda;h),\\
&R_a(\lambda;h)=R_0(\lambda;h)-R_a(\lambda;h) h^2V_a R_0(\lambda;h).
\end{split}\end{equation}
Hence,
\begin{equation}\begin{split}\label{eq:perturb-32}
R(\lambda;h)&-R_0(\lambda;h)-\sum_a (R_a(\lambda;h)-R_0(\lambda;h))\\
&=\sum_a\sum_{b\neq a} R(\lambda;h)h^2
V_b R_a(\lambda;h)h^2 V_a R_0(\lambda;h).
\end{split}\end{equation}
Thus,
the combination of \eqref{eq:perturb-32} and a regularization argument
as in \eqref{eq:fc-tr-reg}, though it is rather cumbersome to keep
track of the regularizing factors,
shows that there is an additional
$h^2$ vanishing here, i.e.\ $a_2=0$.
Note that \eqref{eq:perturb-32} is trace class
when $n=\dim X_0<6$, hence for two-dimensional particles, so in this
case we do not need the regularization procedure.
\end{proof}

We proceed to find the leading coefficient, $d_0=a_4$,
without using the regularization procedure, which
would complicate the calculation. This is justified directly
for two-dimensional particles as mentioned above (the integrand
of \eqref{eq:fc-tr} is trace class); in general, we need to be more careful.
To avoid cumbersome notation, we emphasize the case of $\dim X_0=6$,
i.e.\ the particles are three-dimensional.
First, one can write $R(\lambda;h)$ in a (finite) perturbation series:
\begin{equation}
R(\lambda;h)=R_0(\lambda;h)-R_0(\lambda;h)h^2VR_0(\lambda_0;h)+
R(\lambda;h)h^2VR_0(\lambda;h)h^2VR_0(\lambda_0;h);
\end{equation}
of course the last term can be expanded even further. Similarly,
$R_a(\lambda;h)$ can be expanded in a finite series. The sufficiently
high order terms of the series (e.g.\ from the second term on if
$n=6$, i.e.\ if the particles are three-dimensional),
when $R(\lambda;h)$, etc., are substituted into
the right hand side \eqref{eq:perturb-32}, give trace class terms which
also have at least additional vanishing in $h$ compared to the leading
term, so they do not contribute to $d_0=a_4$. So for $n=6$,
\begin{equation}
\eqref{eq:fc-tr}=
\tr(-\frac{h^4}{2\pi i}
\int\bar\partial_{\lambda}\phit(\lambda)\sum_a\sum_{b\neq a}
R_0(\lambda;h)V_bR_0(\lambda;h)V_aR_0(\lambda;h)
\,d\lambda\wedge d\bar \lambda);
\end{equation}
we a priori know that the integral is trace class since it is the
difference of trace class operators.
In addition, the factors
of $V_a$, $V_b$, can be commuted to the front, and each commutator gives
an extra $h$ vanishing, and lowers the differential order of
the pseudo-differential operator. In fact, the commutators take the
form $[R_0(\lambda;h),h^2V_a]=-R_0(\lambda;h)[h^2\Delta,h^2 V_a]
R_0(\lambda;h)$, so for $n=6$ the corresponding terms
are trace class with extra vanishing
in $h$, so they do not contribute to $d_0$ either. Thus,
\begin{equation}\label{eq:fc-tr-16}
\eqref{eq:fc-tr}=
\tr(-\frac{h^4}{2\pi i}
\int\bar\partial_{\lambda}\phit(\lambda)\sum_a\sum_{b\neq a}
V_bV_aR_0(\lambda;h)^3
\,d\lambda\wedge d\bar \lambda);
\end{equation}
As before, the integrand of \eqref{eq:fc-tr-16} is not trace class,
though the integral is. Now using the Cauchy-Stokes formula, the integral
can be written as
\begin{equation}\label{eq:fc-tr-24}
\eqref{eq:fc-tr}=\tr(\frac{h^4}{2\pi i}
\int \phi(\lambda)(\sum_a\sum_{b\neq a}V_b V_a)
(R_0(\lambda+i0;h)^3-R_0(\lambda-i0;h)^3)\,d\lambda),
\end{equation}
Now the integrand is not trace class, but its trace is well-defined
as a pairing of its kernel with the delta distribution associated to the
diagonal.
We write out this pairing explicitly as the integral of the restriction
of the kernel to the diagonal by writing
$R_0(\lambda+i0;h)^3-R_0(\lambda-i0;h)^3$ as multiplication by a
differentiated
delta distribution
associated to $|\xi|^2=\lambda$, conjugated by the $h$-Fourier transform.
More specifically, we need to calculate
\begin{equation}\label{eq:lead-average-asymp-8}
(2\pi h)^{-n}h^4 (\int_{\Xb_0} (V_a V_b)\,dg)\,((|\xi|^2-(\lambda+i0))^{-3}
-(|\xi|^2-(\lambda-i0))^{-3},1),
\end{equation}
where $(.,1)$ is the distributional pairing of a
compactly supported distribution on $\Rn_\xi$ with $1$. The latter is, in
modified polar coordinates $\xi=\rho^{1/2}\omega$, $|\omega|=1$,
\begin{equation}\label{eq:coeff-16}
\frac{1}{2}\vol(\Sn)((\rho-(\lambda+i0))^{-3}
-(\rho-(\lambda-i0))^{-3},\rho^{(n-2)/2})_\Real.
\end{equation}
Since $(\rho-(\lambda+i0))^{-3}-(\rho-(\lambda-i0))^{-3}=2\pi i\,2^{-1}
\delta_\lambda^{(2)}$ as a distribution in $\rho$, we deduce that
\begin{equation}
\eqref{eq:coeff-16}=\frac{1}{2}\vol(\Sn)2\pi i\,2^{-1}\frac{(n-2)(n-4)}{4}
\lambda^{\frac{n}{2}-3}.
\end{equation}
Hence,
\begin{equation}\label{eq:coeff-32}
\eqref{eq:fc-tr}=
h^{4-n} c(\int_{\Xb_0} (V_a V_b)\,dg)\,
\int\phi(\lambda)\lambda^{\frac{n}{2}-3}\,d\lambda,
\quad c=(2\pi)^{-n}\frac{(n-2)(n-4)}{16}\vol(\Sn).
\end{equation}

In the general, higher dimensional ($n\geq 8$), case one simply has
to keep more terms, one obtains more commutators, and arranges the
factors similarly to how it was done above. Then one still needs
to compute the asymptotics of all the a priori non-trivial terms, but
at this point the only resolvent involved is the free one, and it
is easy to see explicitly that none of the terms except the first one,
namely the one kept in \eqref{eq:fc-tr-24}, contribute to $d_0$,
so the formal calculation that would lead to \eqref{eq:lead-average-asymp-8}
is indeed valid.

So far we have only dealt with the `averaged asymptotics', i.e.\ those
involving $\phi(H/\lambda)$, etc., rather than with the asymptotic behavior
of the $\Cinf$ function $\sigma$ as $\lambda\to+\infty$.
To analyze the latter, recall that
$\|R(\lambda+i0;h)\|_{\bop(H_h^{r,s},H_h^{r+2,-s})}\leq C/h$  for all $r$ and
for all $s>1/2$ on any appropriate weighted
spaces, $\lambda$ in a compact subset of $(0,+\infty)$. Now,
\begin{equation}
R(\lambda)=(H-\lambda)^{-1}=h^2(h^2H-\lambda/h^2)^{-1}=h^2R(\lambda/h^2;h),
\end{equation}
so taking $h=\lambda^{-1/2}$, we deduce that
\begin{equation}\label{eq:high-en-res-est}
\|R(\lambda+i0)\|_{\bop(H_{\lambda^{-1/2}}^{r,s},
H_{\lambda^{-1/2}}^{r+2,-s})}\leq C\lambda^{-1/2},\quad \lambda\geq 1.
\end{equation}
Such estimates follow from
positive commutator estimates (such as the Mourre estimate)
\cite{Jensen:High}, and require
a non-trapping assumption for the semiclassical principal symbol
\cite{Gerard:Semiclassical} -- but this
is just $h^2\Delta$ since the potential is higher order, so the assumption
is automatically satisfied. Microlocal versions of these estimates remain
true. For example, \eqref{eq:WF-B-8} is replaced, for $\supp F\subset (c,c')$
$\supp \tilde F\subset (c',c'')$, $c<c'<c''$, $|c|$, $|c''|$ small, by
\begin{equation}\label{eq:WF-B-8-h}
\|A(F;h) R(\lambda+i0;h) A(\tilde F;h)\|
_{\bop(H_h^{r,s},H_h^{r',s'})}\leq C/h,
\end{equation}
and \eqref{eq:WF-B-16} is replaced, for $\supp\chi\cap\supp\tilde\chi
=\emptyset$, $\supp F\subset (c,c'+\ep)$
$\supp \tilde F\subset (c',c'')$, by
\begin{equation}\label{eq:WF-B-16-h}
\|\chi A(F;h) R(\lambda+i0;h) A(\tilde F;h)\tilde\chi\|
_{\bop(H_h^{r,s},H_h^{r',s'})}\leq C/h.
\end{equation}
In the high-energy version the right hand sides are replaced by
$C\lambda^{-1/2}$. The loss $h^{-1}$ corresponds to the fact that
semiclassically the commutator of two pseudo-differential operators
vanishes to one order higher than the product, see the remarks
following the proof of Proposition~\ref{prop:prop} in the appendix.
Similarly we can apply $R(\lambda+i0)$ iteratively and gain further
decay in $\lambda$, e.g.\ $\|R(\lambda+i0)^2\|\leq C\lambda^{-1}$
on appropriate spaces. As in the preceeding section,
below we apply these results to the resolvent
$(P_\alpha-(\lambda+i0))^{-1}$ of the many-body Hamiltonian $P_\alpha
=\alpha H_L+(1-\alpha H_R)$
on $X_0\times X_0$.

The remark after \eqref{eq:A(h)-ker-norm} shows that
for $s>n/2$, and for all $r\in\Real$,
\begin{equation}\label{eq:psdo-h-est-1}
\|K_{\psi(H/\lambda)}\|_{H^{r,-s}_{\lambda^{-1/2}}}
\leq C\lambda^{n/4},\quad\lambda\geq 1,
\end{equation}
and the same estimate holds for $\delta_{\diag}$ with
$H^{r,-s}_{\lambda^{-1/2}}$ replaced by $H^{-s,-s}_{\lambda^{-1/2}}$.
Similar calculations also yield estimates corresponding to the wave front
set, for example that for $\tilde F$ with $0\nin\supp(1-\tilde F)$,
$A(1-\tilde F)=A(1-\tilde F;\lambda^{-1/2})$ as
in \eqref{eq:q_W-def-h}, and for any $r',s'\in\Real$
\begin{equation}\label{eq:psdo-h-est-2}
\|A(1-\tilde F)K_{\psi(H/\lambda)}\|_{H^{r',s'}_{\lambda^{-1/2}}}
\leq C\lambda^{n/4},\quad\lambda\geq 1,
\end{equation}
and the same estimate holds for $\delta_{\diag}$.
Let
\begin{equation}
\supp\tilde F\subset (c',\infty),\ c'<0,\ |c'| \text{small},
\end{equation}
(we can take,
for example, $c'=-1/2$), $\tilde F$ identically $1$ on a slightly
smaller set, in particular on a neighborhood of $[0,+\infty)$.
Write
\begin{equation}
\Id=A(\tilde F)+A(1-\tilde F)+(\Id-\tilde\psi_0(P_\alpha/\lambda)),
\end{equation}
where $\tilde \psi_0(P_\alpha/\lambda)$ is the operator with amplitude
$p$ in the quantization map \eqref{eq:q_W-def-h}. Since
\begin{equation}\begin{split}\label{eq:semicl-est-32}
&\|(P_\alpha-(\lambda+i0))^{-1}A(\tilde F)K_{\psi(H/\lambda)}\|\leq
\|(P_\alpha-(\lambda+i0))^{-1}A(\tilde F)\|\|K_{\psi(H/\lambda)}\|\\
&\|(P_\alpha-(\lambda+i0))^{-1}A(1-\tilde F)K_{\psi(H/\lambda)}\|\leq
\|(P_\alpha-(\lambda+i0))^{-1}\|\|A(1-\tilde F)K_{\psi(H/\lambda)}\|,\\
&\|(P_\alpha-(\lambda+i0))^{-1}(\Id-\tilde\psi_0(P_\alpha/\lambda))K_{\psi(H/\lambda)}\|\\
&\qquad\qquad
\leq \|(P_\alpha-(\lambda+i0))^{-1}(\Id-\tilde\psi_0(P_\alpha/\lambda))\|
\|K_{\psi(H/\lambda)}\|,
\end{split}\end{equation}
we deduce that for $s>n/2+1$, and for all $r$,
\begin{equation}\label{eq:semicl-est-36}
\|\tilde K_\lambda\|_{H^{r,-s}_{\lambda^{-1/2}}}
=\|(P_\alpha-(\lambda+i0))^{-1}K_{\psi(H/\lambda)}\|_{H^{r,-s}_{\lambda^{-1/2}}}
\leq C\lambda^{n/4-1/2},\quad\lambda>1.
\end{equation}

Now let $F$ be such that $\supp F\subset(-\infty, c')$.
The estimates \eqref{eq:semicl-est-32}, with a factor of $A(F)$ in
front of $(P_\alpha-(\lambda+i0))^{-1}$, yield similar results, but in
$H^{r',s'}_{\lambda^{-1/2}}$ for all $r',s'$, so we deduce that
\begin{equation}
\|A(F)\tilde K_\lambda\|_{H^{r',s'}_{\lambda^{-1/2}}}
\leq C\lambda^{n/4-1/2},\quad\lambda>1.
\end{equation}

Let $b$ be a 2-cluster, and let $\chi=\chi_b$ be a conic cut-off supported
away from $\cup_{a\neq b}\overline{\diag_a}$, but near $\overline{\diag}$.
Let $\tilde\chi$ be
a similar conic cutoff, $\tilde\chi$ identically $1$ on $\supp\chi$.
Then for all $s'$,
\begin{equation}\label{eq:psdo-h-est-3}
\|\tilde\chi(K_{\psi(H/\lambda)}-K_{\psi(H_b/\lambda)})\|
_{H^{r',s'}_{\lambda^{-1/2}}}
\leq C\lambda^{n/4},\quad\lambda\geq 1.
\end{equation}
Moreover, from \eqref{eq:semicl-est-36} and since
$\tilde\chi (\alpha (I_b)_L+(1-\alpha)(I_b)_R)$ is Schwartz,
we deduce that for all $s'$,
\begin{equation}
\|\tilde\chi (\alpha (I_b)_L+(1-\alpha)(I_b)_R) \tilde K_b\|
_{H^{r',s'}_{\lambda^{-1/2}}}
\leq C\lambda^{n/4-1/2},\quad\lambda\geq 1.
\end{equation}
Let $v$ denote the right hand side of \eqref{eq:prop-sing-32}.
Thus,
\begin{equation}
\|\tilde\chi v\|
_{H^{r',s'}_{\lambda^{-1/2}}}
\leq C\lambda^{n/4},\quad
\|A(F)(1-\tilde\chi) v\|
_{H^{r',s'}_{\lambda^{-1/2}}}
\leq C\lambda^{n/4},\quad\lambda\geq 1.
\end{equation}
But
\begin{equation}\begin{split}
\tilde K_\lambda-(\tilde K_b)_\lambda&=(P_\alpha-(\lambda+i0))^{-1}v\\
&=(P_\alpha-(\lambda+i0))^{-1}\tilde\chi v+
(P_\alpha-(\lambda+i0))^{-1}\tilde\psi_0(P_\alpha/\lambda)(1-\tilde\chi)v\\
&\qquad+(P_\alpha-(\lambda+i0))^{-1}A(F)(1-\tilde\chi)v\\
&\qquad+(P_\alpha-(\lambda+i0))^{-1}A(1-F)(1-\tilde\chi)v.
\end{split}\end{equation}
Let $F_0$ satisfy $\supp F_0\subset
(-\infty,\ep)$, $\ep>0$ sufficiently small,
and apply $\chi A(F_0)$ to the previous equation.
We can then estimate each resulting term on the right hand side
in $H^{r',s'}_{\lambda^{-1/2}}$ for all $s'$ by $C\lambda^{n/4-1/2}$.
Namely the first three terms can be estimated by using the operator
estimates of $A(F_0)(P_\alpha-(\lambda+i0))^{-1}$, since it is applied to functions
bounded in $H^{r',s'}_{\lambda^{-1/2}}$ by $C\lambda^{n/4}$ for all $s'$.
On the other hand, the last term
can be estimated by using the operator estimate
of $\chi A(F_0)(P_\alpha-(\lambda+i0))^{-1}A(1-F)(1-\tilde\chi)$. We therefore deduce
that for all $r',s'$,
\begin{equation}\label{eq:psdo-h-est-8}
\|\chi A(F_0) (\tilde K_\lambda-(\tilde K_b)_\lambda))\|
_{H^{r',s'}_{\lambda^{-1/2}}}
\leq C\lambda^{n/4-1/2},\quad\lambda\geq 1.
\end{equation}

Combined with a similar result for $(\tilde K_a)_\lambda-(\tilde K_0)_\lambda$,
and summing over a partition of unity $\chi_b$ of a neighborhood
of $\overline{\diag}$, we deduce a bound for $\chi'A(F_0)\tilde K'_\lambda$
if $\chi'$ is supported near the diagonal. Since the corresponding bound
is automatic away from the diagonal even for $A(F_0)\tilde K_\lambda$
by \eqref{eq:WF-B-16-h},
we conclude that
\begin{equation}\label{eq:prop-sing-h-20}
\|A(F_0)\tilde K'_\lambda\|
_{H^{r',s'}_{\lambda^{-1/2}}}
\leq C\lambda^{n/4-1/2},\quad\lambda\geq 1.
\end{equation}

Our proof of the smoothness of
$\sigma=\langle\tilde K'_\lambda,\delta_{\diag}\rangle$ then first
yields that $|\sigma(\lambda)|\leq C \lambda^{(n-1)/2}$, since
the appropriate microlocal norms of $\tilde K'_\lambda$ are bounded by
$C'\lambda^{n/4-1/2}$, while those of $\delta_{\diag}$ are bounded by
$C''\lambda^{n/4}$.
Indeed, in \eqref{eq:explicit-pairing}, for $F\in\Cinf_c(\Real)$
supported in $(-\infty,\ep)$, identically $1$ near $0$,
and for $s>n/2$ and $\lambda>1$,
\begin{equation}\begin{split}
&|\langle (\Id-\psi_0(P_\alpha/\lambda))
\tilde K'_\lambda,\delta_{\diag}\rangle|
\leq\|(\Id-\psi_0(P_\alpha/\lambda))\tilde K'_\lambda\|_{H^{s,s}}
\|\delta_{\diag}\|_{H^{-s,-s}}\leq C\lambda^{n/2-1/2},\\
&|\langle A(F_0)\tilde K'_\lambda,\delta_{\diag}\rangle|
\leq\|A(F_0)\tilde K'_\lambda\|_{H^{s,s}}\|\delta_{\diag}\|_{H^{-s,-s}}
\leq C\lambda^{n/2-1/2},\\
&|\langle \tilde K'_\lambda,(\psi_0(P_\alpha/\lambda)-A(F_0)^*)\delta_{\diag}\rangle|\\
&\qquad\leq\|\tilde K'_\lambda\|_{H^{-s,-s}}
\|(\psi_0(P_\alpha/\lambda)-A(F_0)^*)\delta_{\diag}\|_{H^{s,s}}
\leq C\lambda^{n/2-1/2}.
\end{split}\end{equation}
The identity
$\lambda dR(\lambda)/d\lambda=HR^2(\lambda)-R(\lambda)$ also allows us
to deduce that
\begin{equation}
|(\lambda\pa_{\lambda})^k\sigma(\lambda)|\leq
C_k \lambda^{(n-1)/2},
\end{equation}
hence $\sigma$ is a symbol (outside a compact set). Applying a perturbation
series argument, as for the `averaged asymptotics' above
with $\lambda$ replaced by $\lambda+i0$,
then shows that for all $k\geq 0$, $\sigma$ has an asymptotic expansion, up to
$\lambda^{(n-k)/2}$, modulo symbols of order $(n-k)/2$.
(Note that a priori there is only a gain of $\lambda^{-1/2}$, rather than
$\lambda^{-1}$, between consecutive terms of the perturbation series,
due to \eqref{eq:high-en-res-est}; this does not appear in the trace
due to the special behavior of the kernel near the diagonal: sharper
(subconic) localization near the diagonal would give, even a priori,
a better result.)
Hence $\sigma$ is indeed a classical symbol, i.e.\ has
an asymptotic expansion, with the top coefficient calculated above.
In view of \eqref{eq:coeff-32},
we have thus proved the theorem from the introduction, which we now
restate.

\begin{thm}\label{thm:3-body}
Suppose that the pair potentials $V_a$ are Schwartz (on $X^a$).
Then the spectral shift
function $\sigma$ is $\Cinf$ on $\Real\setminus\Lambda$.
Moreover,
$\sigma$ is a symbol outside a compact set, and it has a full asymptotic
expansion as $\lambda\to+\infty$:
$$
\sigma(\lambda)\sim \sum_{j=0}^\infty \lambda^{\frac{n}{2}-3-j} c_j,\quad
c_0=C_0\sum_a \sum_{b\neq a}\int_{X_0} V_a V_b\,dg,
$$
where $C_0=\frac{1}{16}(n-2)(n-4)(2\pi)^{-n}\vol(\Sn)$
depends only on $n=\dim X_0$.
\end{thm}

\section{Many-body spectral shift functions}
We now define a modified spectral shift
function in the general $N$-body setting.
Essentially the same proofs as above show its smoothness away from
the thresholds and yield its high energy asymptotics, though now
the combinatorial part becomes a little more complicated. We define
these recursively for subsystems, starting with the free (i.e.\ $N$-)
cluster. So, for $\phi\in\Cinf_c(\Real)$, let
\begin{equation}\begin{split}\label{eq:T-def-8}
&T(X_0,X^0,\calX^0)(\phi)=\phi(H_0),\\
&T(X_0,X^a,\calX^a)(\phi)=\phi(H_a)
-\sum_{c\leq a,\ c\neq a} T(X_0,X^c,\calX^c)(\phi).
\end{split}\end{equation}
The second equation can be rewritten as
\begin{equation}
\phi(H_a)=\sum_{c\leq a} T(X_0,X^c,\calX^c)(\phi),
\end{equation}
and it defines $T(X_0,X^a,\calX^a)(\phi)$ recursively.
Note that $T(X_0,X^a,\calX^a)(\phi)$ is an operator on $X_0$; more
precisely, it is in $\PsiSc^{-\infty,0}(\Xb_0;\calC)$, though some
of the elements of $\calC$ can be dropped.
We write
\begin{equation}\label{eq:T-def-24}
T(\phi)=T(X_0,X^1,\calX^1)(\phi).
\end{equation}
For example, if $H$ is a two-body Hamiltonian, we get
$T(\phi)=\phi(H)-\phi(H_0)$, and if $H$ is a three-body Hamiltonian
we obtain the Buslaev-Merkurev expression, $T(\phi)=\phi(H)-\phi(H_0)
-\sum_{\#a=2}(\phi(H_a)-\phi(H_0))$.
We continue with a lemma.

\begin{lemma}
Suppose that $V_a\in x^k\Cinf(\Xb^a)$.
Then for all $a$, $T(X_0,X^a,\calX^a)(\phi)$
is in $\PsiSc^{-\infty,k}(\Xb_0,\calC)$ away from $C_a$, i.e.\ on
$C_0\setminus C_a$. In particular,
if all $V_a$ are Schwartz, then for all $a$, $T(X_0,X^a,\calX^a)(\phi)$
is in $\PsiSc^{-\infty,\infty}(\Xb_0,\calC)$ away from $C_a$.
\end{lemma}

\begin{proof}
This statement is empty for
$T(X_0,X^0,\calX^0)(\phi)$, since $C_0\setminus C_0=\emptyset$.
We proceed by induction, assuming that we have shown that for all
$c\leq a$, $c\neq a$, $T(X_0,X^c,\calX^c)(\phi)$ is in
$\PsiSc^{-\infty,\infty}(\Xb_0,\calC)$ away from $C_c$.
Suppose that $p\in C_{b,\reg}$, $p\nin C_a$; in particular, $b\not\geq a$
(for then $C_a\supset C_b$ would hold).
Then $T(X_0,X^c,\calX^c)(\phi)$
is in $\PsiSc^{-\infty,\infty}(\Xb_0,\calC)$ near $p$ unless $p\in C_c$,
i.e.\ unless $C_b\subset C_c$, i.e.\ $b\geq c$. Now,
\begin{equation}\begin{split}\label{eq:T-exp-16}
T(X_0,X^a,\calX^a)(\phi)=\phi(H_a)&-\sum_{c\leq a,\ c\neq a\ c\leq b}
T(X_0,X^c,\calX^c)(\phi)\\
&-\sum_{c\leq a,\ c\neq a\ c\not\leq b}T(X_0,X^c,\calX^c)(\phi),
\end{split}\end{equation}
and we have just seen that each term in the last sum is in
$\PsiSc^{-\infty,\infty}(\Xb_0,\calC)$ near $p$.
On the other hand, let $d$ be the maximal element
with the property $d\leq a$ and $d\leq b$. Note that this maximal element
is unique, namely it is given by $X^d=\Span\{X^c:\ c\leq a,\ c\leq b\}\subset
X^a\cap X^b$,
i.e.\ $X_d=\bigcap\{X_c:\ c\leq a,\ c\leq b\}\supset X_a+ X_b$,
which is a collision plane
since $\calX$ is closed under intersections. In particular, $d\neq a$ since
$a\not\leq b$,
hence the first sum in \eqref{eq:T-exp-16} is
$\sum_{c\leq d}T(X_0,X^c,\calX^c)(\phi)=\phi(H_d)$, so
\begin{equation}\label{eq:T-exp-18}
T(X_0,X^a,\calX^a)(\phi)=\phi(H_a)-\phi(H_d)
-\sum_{c\leq a,\ c\neq a,\ c\not\leq b}
T(X_0,X^c,\calX^c)(\phi).
\end{equation}
Since $H_a-H_d=\sum_{c\leq a,\ c\not\leq b}V_c$, it is in $x^k\Cinf(\Xb_0)$
(resp.\ Schwartz) near $p$
if the potentials $V_c$ are in $x^k\Cinf(\Xb^c)$ (resp.\ Schwartz on $X^c$),
hence the local
nature of the construction of $\phi(H_a)$ and $\phi(H_d)$ (functional
calculus and resolvent construction) yields that $\phi(H_a)-\phi(H_d)$
is in $\PsiSc^{-\infty,\infty}(\Xb_0,\calC)$ near $p$, hence providing
the inductive step.
\end{proof}

This shows, in particular, that
$T(\phi)\in\PsiSc^{-\infty,k}(\Xb_0,\calC)$, and is hence trace class if $k>n$.
Moreover, the map $\Cinf_c(\Real)\ni\phi\mapsto
\tr(T(\phi))\in\Cx$ is linear and continuous.
We thus make the following definition.

\begin{Def}
Let $H$ be a many-body Hamiltonian, and define $T$ by
\eqref{eq:T-def-8}-\eqref{eq:T-def-24}. Suppose that the potentials $V_a$
are symbols of order $k>n$ on $X^a$: $V_a\in S^{-k}(X^a)$.
The modified spectral shift function $\sigma$ is defined, as a distribution on
$\Real$, by
\begin{equation}\label{eq:sigma-N-def}
\sigma(\phi)=\tr(T(\phi)),\ \phi\in\Cinf_c(\Real).
\end{equation}
\end{Def}

For Schwartz potentials $V_c$ (on $X^c$),
the arguments presented in the previous sections apply, with the
result that $\sigma$ is $\Cinf$ on
$\Real\setminus\Lambda$. Indeed, note first that the wave front set of
kernel of $T(X_0,X^a,\calX^a)(\phi)$
is in the compressed conormal bundle of the
diagonal (since the operator is a linear combination of the $\phi(H_c)$),
hence $\tau=0$ on it. Let $\psi\in\Cinf_c(\Real)$, and let
$\psit(t)=\psi(t)(t-(\lambda+i0))^{-1}$.
We will write $T(X_0,X^a,\calX^a)(\psit)$ even though $\psit$ is not smooth.
We claim that near any $p\nin C_a$, $T(X_0,X^a,\calX^a)(\psit)$ is a sum of
terms $T_j$, each of which satisfies that for some $a_j\leq a$,
$(H_{a_j}-\lambda)T_j,T_j(H_{a_j}-\lambda)$ are in
$\PsiSc^{-\infty,\infty}(\Xb_0,\calC)$ near $p$. The proof again proceeds
by induction, the statement being empty for $a=0$. So suppose that
for all $c\leq a$, $c\neq a$, we have shown the claim, and
suppose that $p\in C_{b,\reg}$. Now
\eqref{eq:T-exp-18} becomes
\begin{equation}\label{eq:T-exp-24}
T(X_0,X^a,\calX^a)(\psit)=\psit(H_a)-\psit(H_d)
-\sum_{c\leq a,\ c\neq a,\ c\not\leq b}
T(X_0,X^c,\calX^c)(\psit).
\end{equation}
Each term in the sum on the right hand side can be written as a sum
of operators $T_{c,j}$ with the desired properties by the induction
hypothesis. On the other hand,
\begin{equation}
(H_a-\lambda)(\psit(H_a)-\psit(H_d))
=\psi(H_a)-\psi(H_d)+(\sum_{c\leq a,\ c\not\leq d}V_c)\psi(H_d).
\end{equation}
Since $c\leq b$, together with $c\leq a$, would imply $c\leq d$, we
deduce that each $V_c$ is Schwartz near $p$, hence the last term
is in $\PsiSc^{-\infty,\infty}(\Xb_0,\calC)$ near $p$, while the
same statement for $\psi(H_a)-\psi(H_d)$ has already been demonstrated
above. Proceeding in the three-body setting, i.e.\ using the
propagation of singularities, shows that the wave front
set of the kernel of $T(\psit)$ is in $\tau\leq \tau_0<0$, and then we
deduce that $\sigma$ is continuous. Iterating the argument gives that $\sigma$
is $\Cinf$.

The pseudodifferential functional calculus shows, as in the previous section,
that for $\phi\in\Cinf_c(\lambda)$,
\begin{equation}
\sigma(\phi(./\lambda))\sim\sum_{j=0}^\infty \lambda^{(n-j)/2} a_j,
\end{equation}
with $a_j=0$ for all odd $j$.
The symbol estimates for the $\Cinf$ function $\sigma$ itself also proceed
as before, so $|(\lambda\pa_{\lambda})^k\sigma(\lambda)|
\leq C_k\lambda^{(n-1)/2}$. A perturbation series argument again yields
a full asymptotic expansion. Since $R_a(\lambda)-R_0(\lambda)=R_a(\lambda)
(\sum_{c\leq a}V_c)R_0(\lambda)$, $\lambda\nin\Real$, we deduce that
the leading terms $a_0$ and $a_1$ vanish.
If all potentials are pair potentials,  i.e.\ if
$V_c=0$ for all clusters $c$ that are not $N-1$-clusters, one can show that
$a_j =0$ for $j\le 2N-1$ and $a_{2N }= c_{n,N}C(V)$
where
$C(V) = \sum_{\sigma} \int V_{\sigma_1} \cdots V_{\sigma_{N-1}} dx$,
the sum is taken over all sequences $\sigma = (\sigma_1, \sigma_2,
\cdots, \sigma_{N-1})$ of $(N-1)$ -cluster decompositions $\sigma_j$
such that  $\cup_{j=1}^{N-1} \sigma_j = a_{\max}$
and $c_{n,N}$ is a contant depending only on $n$ and $N$
which can be calculated by the method of Section~\ref{sec:high}.

We have thus proved the following result.

\begin{thm}\label{thm:N-body-trace}
Suppose that for all $a$, $V_a$ is a Schwartz function on $X^a$.
The modified spectral shift function, $\sigma$, defined as a distribution by
\eqref{eq:sigma-N-def}, is $\Cinf$ on $\Real\setminus\Lambda$. In addition,
$\sigma$ is a symbol outside a compact set, and it has a full asymptotic
expansion:
$$
\sigma(\lambda)\sim\sum_{j=0}^\infty\lambda^{\frac{n}{2}-2-j}c_j,
\quad\lambda\to+\infty.
$$
Here we do not assume that all potentials are pair potentials.
\end{thm}

\appendix
\section{Sketch of relevant positive commutator estimates}

To illustrate the proof of the propagation of singularities by
positive commutator estimates in \cite{Vasy:Bound-States},
we sketch the proof of a simpler version
here which still suffices for the purposes of the present paper.
We state it as a wave front set estimate, though as we see below,
as usual with positive commutator estimates, it
actually amounts to a microlocal energy estimate.

First we introduce some notation. For an cluster $a$ (possibly $a=0$!)
let
\begin{equation}
x_a=|w_a|^{-1},\ y_a=w_a/|w_a|,\ z_a=w^a/|w_a|;
\end{equation}
these are local coordinates on the radial compactification $\bar X_0$ near
$\bar X_a$. Let
\begin{equation}
\tau=-\frac{w\cdot\xi}{|w|},\ \tau_a=-\frac{w_a\cdot\xi_a}{|w_a|},
\end{equation}
so $\tau=\tau_a$ at $\bar X_a$. We also write $x=|w|^{-1}$. In addition,
it is sometimes convenient, for $X_b\supset X_a$, i.e.\ $X^b\subset X^a$,
to decompose $w^a\in X^a$ as $(w^{ab},w^b_a)\in X^b\oplus (X^a\ominus X^b)$,
and write $(z_{ab},z^b_a)$ accordingly. We also let $\scHg^c$ be the
rescaled Hamilton vector field of $\Delta_{X_c}$, so
\begin{equation}
\scHg^c=2\frac{\xi_c}{|w_c|}\cdot\pa_{w_c};
\end{equation}
this should be regarded as a vector field on $T^*X_c$ which extends
to a smooth vector field on $\sct\Xb_c$ tangent to the boundary
$\sct_{C_c}\Xb_c$, and hence can be considered as a vector field on
$\sct_{C_c}\Xb_c$.
Recall also that the part of the characteristic variety corresponding
to the bound states of $H^b$ is
\begin{equation}
\Sigma_b(\lambda)=\{\zeta=(y_b,\xi_b)\in\sct_{C_b}\Xb_b:
\ \lambda-|\xi_b|^2\in\pspec{H^b}\}\subset\sct_{C_b}\Xb_b.
\end{equation}
For $C_a\subset C_b$, we write the projection
\begin{equation}
\hat\pi_{ba}:\Sigma_b(\lambda)\cap\sct_{C_{a,\reg}}\Xb_b
\to\sct_{C_{a,\reg}}\Xb_a\subset\dot\Sigma(\lambda);
\end{equation}
this is the restriction of
$\pi_{ba}:\sct_{C_{a,\reg}}\Xb_b
\to\sct_{C_{a,\reg}}\Xb_a$ to $\Sigma_b(\lambda)$.
For $A\in\PsiSc^{-\infty,l}(\Xb,\calC)$, the operator wave front set
$\WFScp(A)$ was defined in \cite{Vasy:Propagation-Many}
as a subset of the compressed cotangent bundle
$\scdt \Xb$.
Namely, let $p$ denote the projection
$[\Xb;\calC]\times\Xb^*\to\scdt\Xb$ given by the composition of the blow-down
map and the projection $\pi:\sct\Xb\to\scdt\Xb$. Then
$\zeta\in\scdt\Xb\setminus\WFScp(A)$ means that in a neighborhood
of $p^{-1}(\{\zeta\})$, the amplitude $a$
defining $A$ as in \eqref{eq:cl-symb-def} vanishes to infinite order.
This notion is independent of the choice of quantization.

The main technical result is thus the following.

\begin{prop}\label{prop:prop}
\cite[Weaker version of Proposition 7.1]{Vasy:Bound-States}
Suppose that $H$ is a many-body Hamiltonian.
Let $u\in\dist(\Xb)$,
$\lambda\nin\Lambda_1$.
Let $\bar y_a\in C_{a,\reg}$, $\bar\tau$ such that
$\lambda-\bar\tau^2\nin\Lambda_1$.
Suppose that for all $\bar\zeta=(\bar y_a,\bar \xi_a)\in\sct_{C_a}\Xb_a$ with
$\tau(\bar\zeta)=\bar\tau$, we have
$\bar\zeta\nin\WFSc((H-\lambda)u)$.
Then there exist $\bar\delta_0>0$ and $C>0$, depending only on $\bar\zeta$,
with the following property.
For all $\delta_0\in(0,\bar\delta_0)$
such that if
\begin{equation}\label{eq:prop-ii}
\forall\zeta\ \text{s.t.}\ |y_0(\zeta)-\bar y_a|<C\delta_0
\Mand \bar\tau+\delta_0/3<\tau(\zeta)<\bar\tau+\delta_0
\Rightarrow\zeta\nin\WFSc(u),
\end{equation}
then $\bar\zeta\nin\WFSc(u)$.

In fact, there exists $C'>0$ so that
the following holds. For any $r,s\in\Real$, $r'>-r$, $s'>-s$,
depending only on $\bar\zeta$, there exists $C_1>0$ such that
for all $u$ as above,
\begin{equation}\begin{split}\label{eq:prop-exp-est-8}
\|A(\chi_{U'}(y_0)&\chi_{(-\delta_0/3,\delta_0/3)}(\tau-\bar\tau))u\|_{H^{r,s}}
\leq C_1(\|A(\chi_U(y_0)\chi_{(\delta_0/3,\delta_0)}(\tau-\bar\tau))u\|_{H^{r,s}}\\
&+\|A(\chi_U(y_0)\chi_{(-\delta_0/3,\delta_0)}(\tau-\bar\tau))
(H-\lambda)u\|_{H^{r,s+1}}
+\|u\|_{H^{r',s'}}),
\end{split}\end{equation}
where $U$ is the ball $|y_0(\zeta)-\bar y_a|<C\delta_0$,
$U'$ the ball $|y_0(\zeta)-\bar y_a|<C'\delta_0$, $\chi_U(y)$,
resp.\ $\chi_I(\tau-\bar\tau)$ denote
a smoothed characteristic function of $U_y$, resp.\ the interval
$I_{\tau-\bar\tau}$, and $A$ denotes quantization as in \eqref{eq:q_W-def}.
\end{prop}

\begin{rem}
Since $\bar\zeta\nin\WFSc((H-\lambda)u)$, by elliptic regularity
we deduce that $\bar\zeta\nin \WFSc(u)$
for $\bar\zeta\nin\dot\Sigma(\lambda)$, i.e.\ for $\bar\zeta$ not in the
$\lambda$ characteristic set (energy shell).

The estimate of \eqref{eq:prop-exp-est-8} implies \eqref{eq:WF-B-16}
directly. Indeed, if $\tilde\chi$ is such that $\supp\chi_{U'}\cap\supp
\tilde\chi=\emptyset$, consider $u=R(\lambda+i0)v$,
$v=A(\tilde\chi(y_0)
\chi_{(-\delta_0/4,\delta_0/4)}(\tau-\bar\tau))f$, $f\in H^{r',s'}$. Then 
the second and third terms are directly bounded in terms of $v$,
while the first term is bounded in terms of $f$ due to the
boundedness of $A(\chi_{(\delta_0/3,\delta_0)}(\tau-\bar\tau))
R(\lambda+i0)A(\chi_{(-delta_0/4,\delta_0/4)}(\tau-\bar\tau))$
between any two weighted Sobolev spaces.
\end{rem}

\begin{proof}(Sketch, see
\cite[Proof of Proposition 7.1]{Vasy:Bound-States}
for complete details.)
We give the full commutator construction at the symbol level,
and indicate why it gives rise to a microlocally positive commutator.
In fact, the commutator will be positive in part of phase space, negative
(or not necessarily positive) in another part of phase space. The
propagation of singularities estimates, which should be thought of
as microlocal energy estimates, work by estimating $u$ in the former
region in terms of $u$ in the latter region and $(H-\lambda)u$ in the
union of both regions.

Employing an iterative argument, we may assume
that for all $\bar\zeta=(\bar y_a,\bar\xi_a)$ with $\tau(\bar\zeta)=\bar\tau$,
$\bar\zeta\nin\WFSc^{*,l}(u)$, and we need
to show that $\bar\zeta\nin\WFSc^{*,l+1/2}(u)$. (We can start the induction
with an $l$ such that $u\in H^{*,l}(X)$.)

For points $\zeta=(y_b,\xi_b)$ in $\sct_{C_b}\Xb_b$,
\begin{equation}
-\scHg^b\tau(\zeta)=2(|\xi_b|^2-\tau(\zeta)^2)\geq 0.
\end{equation}
For $\zeta\in \Sigma_b(\lambda)\subset\sct_{C_b}\Xb_b$,
$\lambda-|\xi_b|^2=\ep_\beta\in\pspec(H^b)$, so
\begin{equation}\label{eq:prop-16}
-\scHg^b\tau(\zeta)=2(\lambda-\tau(\zeta)^2-\ep_\beta).
\end{equation}
We define
\begin{equation}
c_0=\frac{1}{2}
\inf\{-\scHg^b\tau(\bar\zeta_b):\exists\bar\xi_a\Mst
\tau((\bar y_a,\bar\xi_a))=\bar\tau,\ \bar\zeta_b\in
\pih_{ba}^{-1}((\bar y_a,\bar\xi_a)),
\ C_b\supset C_a\}.
\end{equation}
Due to \eqref{eq:prop-16}, and due to $\lambda-\bar\tau^2\nin\Lambda_1$,
we deduce that $c_0>0$.
Thus, there exists $\delta_1>0$ such that for all clusters $b$ with
$C_b\supset C_a$, and for all $\zeta\in\Sigma_b(\lambda)$
that satisfies $|y_a(\zeta)-\bar y_a|<\delta_1$,
$|z_a(\zeta)|<\delta_1$,
$|\tau(\zeta)-\bar\tau|<\delta_1$, we deduce that
\begin{equation}\label{eq:pos-comm-loc-32}
-\scHg^b\tau(\zeta)\geq 3c_0/2.
\end{equation}

Our positive commutator estimates will arise by considering functions
\begin{equation}
\phi=\bar\tau-\tau+\frac{\beta}{\ep}(|z_a|^2+|y_a-\bar y_a|^2),
\end{equation}
where $\beta>0$ will be fixed later and $\ep>0$ is arbitrary as long as it
is sufficiently small.
Note that for all $b$ with $C_b\supset C_a$,
$\scHg^b|z_a|^2=4z_{ab}\cdot\xi^a_b$ under the decomposition $\xi_b=(\xi_a,
\xi^a_b)$, so $\scHg^b|z_a|^2\leq C_1|z_a|$ on $\Sigma_b(\lambda)$,
and similarly, possibly by increasing $C_1$,
$|\scHg^b(y_a-\bar y_a)^2|\leq C_1|y_a-\bar y_a|$.

Now suppose that
\begin{equation}\label{eq:case-ii-cond-1}
\phi\leq 2\ep,\ \bar\tau-\tau\geq -2\ep.
\end{equation}
Then we
conclude that
\begin{equation}\label{eq:supp-q-a-2}
|\bar\tau-\tau|\leq 2\ep,\ |z_a|\leq 2\ep/\sqrt{\beta},
\ |y_a-\bar y_a|\leq 2\ep/\sqrt{\beta}.
\end{equation}
Let $\beta=(c_0/4C_1)^2$. For $\ep>0$ small, \eqref{eq:case-ii-cond-1}
thus implies that $|y_a(\zeta)-\bar y_a|<\delta_1$,
$|z_a(\zeta)|<\delta_1$,
$|\tau(\zeta)-\bar\tau|<\delta_1$, so we deduce from
\eqref{eq:pos-comm-loc-32} that
\begin{equation}
\scHg^b\phi\geq -\scHg^b\tau-2\sqrt{\beta} C_1
>c_0,\ \text{where}\ \beta=(c_0/4C_1)^2.
\end{equation}

The positive commutator estimate then arises by considering the following
symbol $q$ and quantizing it as in \eqref{eq:q_W-def}.
Let $\chi_0\in\Cinf(\Real)$ be equal to $0$ on $(-\infty,0]$ and
$\chi_0(t)=\exp(-1/t)$ for $t>0$. Thus, $\chi_0'(t)=t^{-2}\chi_0(t)$, $t>0$,
and $\chi'_0(t)=0$, $t\leq 0$.
Let $\chi_1\in\Cinf(\Real)$ be $0$
on $(-\infty,0]$, $1$ on $[1,\infty)$, with $\chi_1'\geq 0$
and $\chi_1(t)=\exp(-1/t)$
on some small interval
$(0,t_0)$, $t_0>0$.
Furthermore, for $A_0>0$ large, to be determined, let
\begin{equation}\label{eq:prop-22b}
q=\chi_0(A_0^{-1}(2-\phi/\ep))\chi_1((\bar\tau-\tau)/
\ep+2).
\end{equation}
Thus, $q(\tilde\zeta)=\chi_0(2/A_0)>0$, and on $\supp q$ we have
\begin{equation}\label{eq:supp-q-a}
\phi\leq 2\ep\Mand \bar\tau-\tau\geq -2\ep,
\end{equation}
which is \eqref{eq:case-ii-cond-1},
so $\supp q$ is a subset of \eqref{eq:supp-q-a-2}.
We also see that
as $\ep$ decreases, so does $\supp q=\supp q_\ep$,
in fact, if $0<\ep'<\ep$ then $q_\ep>0$ on $\supp q_{\ep'}$.
Note that
by reducing $\ep$,
we can make $q$ supported in an arbitrary small
neighborhood of $\bar y_a$ and $\bar\tau_a$.

Let $\psit\in\Cinf_c(\Real)$ be identically $1$ near $0$ and supported
close to $0$. We also define
\begin{equation}
\tilde q=\psit(x)q.
\end{equation}
Let $A$ be the operator given by \eqref{eq:q_W-def}
with $\tilde q$ in place of $q$. Note that this includes a spectral
cutoff in the definition of $A$.

The commutator $i[\Delta_{X_c},A]$ is given to top order by $\scHg^c q$.
This is the commutator that gives microlocal positivity on the
$L^2$ eigenspace of $H^c$, see e.g.\ the Froese-Herbst proof of
the Mourre estimate \cite{FroMourre}. We proceed to esimate $\scHg^c q$
directly.

Thus,
\begin{equation}\begin{split}\label{eq:scHg q-calc}
\scHg^c q=-A_0^{-1}&\ep^{-1}\chi'_0(A_0^{-1}(2-\phi/\ep))
\chi_1((\bar\tau-\tau)/\ep+2)\scHg^c \phi\\
&-\ep^{-1}
\chi_0(A_0^{-1}(2-\phi/\ep))\chi_1'((\bar\tau-\tau)/\ep+2)
\scHg^c \tau.
\end{split}\end{equation}
Then
\begin{equation}\label{eq:scHg q-pos}
\scHg^c q=-\tilde b^2_c+e_c
\end{equation}
with
\begin{equation}
\tilde b^2_c=A_0^{-1}\ep^{-1}\chi'_0(A_0^{-1}(2-\phi/\ep))
\chi_1((\bar\tau-\tau)/\ep+2)\scHg^c \phi.
\end{equation}
Hence, with
\begin{equation}
b^2=c_0 A_0^{-1}\ep^{-1}\chi'_0(A_0^{-1}(2-\phi/\ep))
\chi_1((\bar\tau-\tau)/\ep+2),
\end{equation}
we deduce that
\begin{equation}\label{eq:prop-32}
\scHg^c q\leq -b^2+e_c.
\end{equation}
Moreover,
\begin{equation}\label{eq:prop-33}
b^2\geq (c_0 A_0/16) q
\end{equation}
since $\phi\geq \bar\tau-\tau\geq -2\ep$ on $\supp q$, so
\begin{equation}\begin{split}\label{eq:prop-34}
\chi'_0(A_0^{-1}(2-\phi/\ep))
&=A_0^2(2-\phi/\ep)^{-2}\chi_0(A_0^{-1}(2-\phi/\ep))\\
&\geq (A_0^2/16)
\chi_0(A_0^{-1}(2-\phi/\ep)).
\end{split}\end{equation}
On the other hand, $e_c$ is supported where
\begin{equation}\label{eq:supp-e-fine}
-2\epsilon\leq\bar\tau-\tau\leq-\ep,
\ |y_a-\bar y_a|,\ |z_a|\leq 2\ep/\sqrt{\beta}.
\end{equation}
By our assumption, this region is disjoint from $\WFSc(u)$, if we
choose $\ep>0$ sufficiently small.
Moreover, by \eqref{eq:supp-q-a},
for $\ep>0$ sufficiently small, we deduce from the
inductive hypothesis that $\supp q$ (hence $\supp b$) is disjoint
from $\WFSc^{*,l}(u)\cap\dot\Sigma(\lambda)$.

Let $B\in\PsiSc^{-\infty,0}(\Xb,\calC)$ be a quantization of
$bq^{1/2}$ as in \eqref{eq:q_W-def}. Suppose
that $M>0$ and $\ep'>0$. By choosing $A_0$ large, depending on $M$, $\ep'$,
(using \eqref{eq:prop-34}), one can derive a positive commutator estimate from
\eqref{eq:prop-32} using the many-body pseudo-differential calculus,
see \cite[End of proof of Proposition~7.1]{Vasy:Bound-States} for details.
Apart from technical details it essentially corresponds to using the Mourre
estimate and the functional calculus microlocally, namely that when
localized in phase space in the region of interest, the commutator of
a quantization of $\phi$ is positive. We deduce that
there exists $\delta'>0$, such that for $\psi\in\Cinf_c(\Real)$ is supported in
$(\lambda-\delta',\lambda+\delta')$, $\psi\equiv 1$ near $\lambda$,
$E\in\PsiSc^{-\infty,0}(\Xb,\calC)$, $F\in\PsiSc^{-\infty,1}(\Xb,\calC)$
with $\WFScp(E),\WFScp(F)$ in a small neighborhood
of $\dot\Sigma(\lambda)$,
\begin{equation}\label{eq:comm-51a}
\WFScp(E)\subset\supp e,\ \WFScp(F)\subset\supp q,
\end{equation}
such that
\begin{equation}\label{eq:comm-53a}
i\psi(H)x^{-1/2}[A^*A,H]x^{-1/2}\psi(H)-M\psi(H)A^*A\psi(H)
\geq (2-2\ep')\psi(H)B^*B\psi(H)+E+F.
\end{equation}
By $\supp e$ we mean the support of the function
$\chi_0(A_0^{-1}(2-\phi/\ep))\chi_1'((\bar\tau-\tau)/\ep+2)$, which is
independent of $c$ in \eqref{eq:scHg q-calc}.
Here $F$ is the error term, it has first order decay, hence it is
`negligible'. On the other hand, $E$ has the same order as $B^*B$, and
it is negative (i.e.\ has the opposite sign of $B^*B$) in part of the
phase space. As mentioned above,
positive commutator estimates for approximate
solutions $u$, i.e.\ $(H-\lambda)u$ microlocally Schwartz, work by
estimating
$\|B\psi(H)u\|^2$ in terms of $\langle u,Eu\rangle$ (plus error terms),
i.e.\ $u$ is estimated on $\supp b$ by its estimate on $\supp e$.

One can now use $M$, chosen sufficiently large, to deal with arbitrary
weights $x^{-l-1/2}$. A standard commutator and
regularization argument then
proves that $x^{-l-1/2}Bu\in L^2(X_0)$, which in turn finishes the proof.
We refer to \cite[Proposition~7.1]{Vasy:Bound-States} for details.

Instead of following this route, we prove the corresponding resolvent
estimate. So suppose that
\begin{equation}
u^+_t=(H-(\lambda+it))^{-1}f,\ t>0,
\end{equation}
and $\WFSc(f)$ is disjoint from the region of interest, and
it is, say, in $\tau<\tau_0$, $\tau_0>0$ sufficiently small, so
that $u^+_t$ converges to $(H-(\lambda+i0))^{-1}f$ as $t\to 0$ in
sufficiently large weighted Sobolev spaces. As above,
assume that $u_t^+$ is uniformly bounded in the region of interest in
$H^{*,l}(X_0)$; we want to prove that it is also uniformly bounded
in $H^{*,l+1/2}(X_0)$.
For $\psi\in\Cinf_c(\Real;[0,1])$ supported sufficiently close to
$\lambda$, with $A_l=A\psi(H)x^{-l-1}$, $B_l=x^{-l-1/2}B\psi(H)$,
we deduce from
\eqref{eq:comm-53a} that
\begin{equation}\label{eq:GIS-7}
ix^{l+1/2}[A^*_l A_l,H]x^{l+1/2}\geq x^{l+1/2}((2-2\ep')B_l^*B_l
+E_l+F_l)
x^{l+1/2},\quad\ep'>0,
\end{equation}
$E_l\in\PsiSc^{-\infty,-2l-1}(\Xb,\calC)$,
$F_l\in\PsiSc^{-\infty,-2l}(\Xb,\calC)$, with similar properties
as in \eqref{eq:comm-51a}.
Since
\begin{equation}
\langle u^+_t,i[A_l^*A_l,H]u^+_t\rangle
=-2\im\langle u^+_t,A^*_lA_l(H-(\lambda+it))u^+_t
\rangle-2t\|A_l u^+_t\|^2,
\end{equation}
we conclude that
\begin{equation}\label{eq:GIS-17}
\|B_l u^+_t\|^2+2t\|A_l u^+_t\|^2\leq
|\langle u^+_t,E_l u^+_t\rangle|+|\langle u^+_t,F_l u^+_t\rangle|
+2|\langle u^+_t,A_l^*A_l(H-(\lambda+it))u^+_t\rangle|.
\end{equation}
Since $t>0$, the second term on the left hand side can be dropped.
Since $u^+_t\to u_+$ in $H^{0,l'}(X_0)$ for $l'<-1/2$, we conclude
that for $l\in(-1,-1/2)$ the right hand side stays bounded as $t\to 0$,
for $u^+_t$ is uniformly bounded in $H^{0,l+1/2}(X_0)$ on $\WFScp(E_l)$ and
it is uniformly bounded in $H^{0,l}(X_0)$ on $\WFScp(F_l)$.
Thus, $B_l u^+_t$ is uniformly bounded in $L^2(X_0)$, and as $u^+_t\to
u_+$ in $H^{0,l'}(X_0)$, we conclude that $B_lu_+\in L^2(X_0)$.

\end{proof}

The semiclassical version of the estimate \eqref{eq:prop-exp-est-8} is
\begin{equation}\begin{split}\label{eq:prop-exp-est-h-8}
\|A(\chi_{U'}(y_0)&\chi_{(-\delta_0/3,\delta_0/3)}(\tau-\bar\tau))u\|
_{H^{r,s}_h}
\leq C_1(\|A(\chi_U(y_0)\chi_{(\delta_0/3,\delta_0)}(\tau-\bar\tau))u\|
_{H^{r,s}_h}\\
&+h^{-1}\|A(\chi_U(y_0)\chi_{(-\delta_0/3,\delta_0)}(\tau-\bar\tau))
(H-\lambda)u\|_{H^{r,s+1}_h}
+\|u\|_{H^{r',s'}_h}).
\end{split}\end{equation}
The proof of this proceeds just as above. Equation \eqref{eq:comm-53a}
is replaced by
\begin{equation}\begin{split}
&i\psi(H)x^{-1/2}[A^*A,H(h)]x^{-1/2}\psi(H)-hM\psi(H(h))A^*A\psi(H(h))\\
&\qquad\geq (2-2\ep')h\psi(H(h))B^*B\psi(H(h))+hE+h^2 F,
\end{split}\end{equation}
i.e.\ the principal terms $B^*B$ and $E$ have an extra factor of $h$
(since the semiclassical calculus is commutative to top order in $h$),
and the error term $F$ has a gain of $h$. Note that
$H(h)=\Delta+h^2V$ shows that semiclassically $V$ is two orders lower
in $h$ that $\Delta$, which in fact significantly simplifies the
argument that turns \eqref{eq:prop-32} into \eqref{eq:comm-53a} (for
sufficiently small $h$). Then \eqref{eq:GIS-17} becomes, after dropping
the second term on the left hand side and multiplying through by $h^{-1}$,
\begin{equation}
\|B_l u^+_t\|^2\leq
|\langle u^+_t,E_l u^+_t\rangle|+h|\langle u^+_t,F_l u^+_t\rangle|
+2h^{-1}|\langle u^+_t,A_l^*A_l(H-(\lambda+it))u^+_t\rangle|,
\end{equation}
and then one can finish the proof as before, using that $u^+_t$ is bounded
by $Ch^{-1}$.

\bibliographystyle{plain}
\bibliography{sm}

\end{document}